\documentclass[11pt]{article}
\usepackage[T1]{fontenc}
\usepackage[utf8]{inputenc}
\usepackage{lmodern}
\usepackage{microtype}
\usepackage[a4paper,margin=2.25cm]{geometry}
\usepackage{amsmath,amssymb,amsthm,mathtools}
\usepackage{booktabs,array,longtable}
\usepackage{graphicx}
\usepackage{tikz}
\usetikzlibrary{arrows.meta,positioning,calc,decorations.pathreplacing,fit}
\usepackage{pgfplots}
\pgfplotsset{compat=1.18}
\usepackage{xcolor}
\usepackage{enumitem}
\usepackage{hyperref}
\usepackage[nameinlink,noabbrev]{cleveref}
\usepackage{setspace}
\usepackage{caption}
\usepackage{subcaption}
\usepackage{bm}
\usepackage{mathrsfs}

\definecolor{deepblue}{RGB}{17,72,120}
\definecolor{teal}{RGB}{0,112,112}
\definecolor{burgundy}{RGB}{125,32,48}
\definecolor{gold}{RGB}{160,112,0}
\definecolor{softgray}{RGB}{245,247,249}
\definecolor{greenx}{RGB}{35,120,75}
\hypersetup{
  colorlinks=true,
  linkcolor=deepblue,
  citecolor=burgundy,
  urlcolor=teal,
  pdfauthor={Tony Newton and Juan Antonio Vega Coso},
  pdftitle={Stochastic Orders of the Exit Pair and the Monotone Response to Resetting}
}
\setlist[itemize]{topsep=0.25em,itemsep=0.2em}
\setlist[enumerate]{topsep=0.25em,itemsep=0.2em}
\newtheorem{theorem}{Theorem}[section]
\newtheorem{proposition}[theorem]{Proposition}
\newtheorem{lemma}[theorem]{Lemma}
\newtheorem{corollary}[theorem]{Corollary}
\newtheorem{definition}[theorem]{Definition}
\newcommand{\TP}{\mathrm{TP}_2}
\theoremstyle{remark}
\newtheorem{remark}[theorem]{Remark}

\newcommand{\Pp}{\mathbb P}
\newcommand{\Ee}{\mathbb E}

\newcommand{\ind}{\mathbf 1}

\newcommand{\Cov}{\operatorname{Cov}}

\newcommand{\sgn}{\operatorname{sgn}}
\newcommand{\dd}{\,\mathrm d}

\newcommand{\lr}{\mathrm{lr}}
\newcommand{\rh}{\mathrm{rh}}
\newcommand{\st}{\mathrm{st}}
\newcommand{\pgf}{\mathrm{pgf}}
\newcommand{\Law}{\mathcal L}

\begin{document}

\begin{center}
{\LARGE\bfseries Stochastic Orders of the Exit Pair\\[2mm]
and the Monotone Response to Resetting\par}
\vspace{0.8em}
\begin{tabular}{c@{\hspace{2.2cm}}c}
\textbf{Tony Newton} & \textbf{Juan Antonio Vega Coso}\\
Newton Astro Labs\textsuperscript{*} & IUFFyM, Universidad de Salamanca\\
London, UK & Salamanca, Spain\\
\href{mailto:tony.newton79@gmail.com}{tony.newton79@gmail.com} & \href{mailto:jantovc@gmail.com}{jantovc@gmail.com}\\
\end{tabular}
\vspace{0.45em}

{\footnotesize \textsuperscript{*}Newton Astro Labs is the trading name under which Tony Newton conducts independent computational research in the United Kingdom; it is not a limited company.}

\vspace{0.4em}
\end{center}
\begin{abstract}
Resetting a process, that is, restarting it from its initial position at random times, changes the probability that it leaves an interval through one end rather than the other. We show that the direction of this change is governed by how the exit time depends on the exit side, and we express this through classical stochastic orders. For random walks and birth--death chains we determine when resetting more intensively always moves the exit probability in the same direction, and we show that this can fail: we give an explicit chain in which the direction reverses twice over time, and bound the number of such reversals spectrally. We also show how the direction depends on the starting point. In the continuum limit we obtain the corresponding results for diffusions with drift.
\end{abstract}

\vspace{0.5em}

\section{Introduction and main results}

Stochastic resetting interrupts a process and restarts it from a prescribed state.  Since the foundational diffusion model of Evans and Majumdar \cite{EvansMajumdar2011}, resetting has been used to control mean first-passage times, stationary laws, search processes and absorption probabilities.  In a finite interval, Pal and Prasad derived detailed first-passage and splitting formulas for Brownian motion with Poisson resetting \cite{PalPrasad2019}; renewal formulations and more general reset laws are now a standard part of the subject, and the effect of restart on the probability of success of a Bernoulli trial was analysed by Belan \cite{Belan2018}; see \cite{EvansMajumdarSchehr2020} for a review.  Closest to the present work, Vidmar characterised the laws of a single completion time that are left unchanged, or uniformly ordered, by every restart law, for arbitrary, deterministic and constant-rate restart \cite{Vidmar2022}; in particular, arbitrary and deterministic restart lead to the same characterisations, whereas constant-rate restart does not.  For processes with several possible outcomes, Pal, Dagdug, Ghosh, Boyer and Pal gave a criterion, in terms of the means and coefficients of variation of the conditional completion times, for whether resetting favours a selected outcome \cite{PalEtAl2025}.  The present paper differs from both in its object and in the type of statement.  The object is the pair of conditional laws of a two-outcome first passage, not a single completion time, and the observable is the splitting probability.  The statements are global in the reset mechanism: each classical stochastic order between the two conditional laws is matched with a monotone response over a natural class of reset clocks.  A division similar to Vidmar's appears here: a universal sign under deterministic deadlines alone already forces the usual stochastic order, as it does under arbitrary clocks, while geometric resetting is governed by a different, generating-function order.  Related discrete models include the semi-deterministic reset walk of Villarroel, Montero and Vega \cite{VillarroelMonteroVega2021}, and the midpoint invariance of the ruin probability under resetting is studied in \cite{Vega2026Critical}.

The earlier spectral programme on critical invariants, reset-neutral distributions, resolvent intertwining, and spectral duality provides a complementary perspective \cite{Vega2026Critical,Vega2026Duality,Vega2026Resolvent,Vega2026Algebra}.  Here the main object is instead the joint exit-time/exit-side law and its order-theoretic interaction with independent reset clocks.

This paper addresses the following question: \emph{if two reset mechanisms differ only in how strongly they weight early versus late absorption, when must the ruin probability move monotonically, and in which direction?}

For a nearest-neighbour walk on $\{0,1,\ldots,a\}$ with absorbing endpoints, reset to its initial site $z$ after an independent renewal time, the eventual lower-boundary probability is a weighted average of the exact-time conditional exit probabilities $\varphi_n=\Pp_z(J=0\mid T=n)$ of the reset-free walk, where $T$ and $J$ are the absorption time and the exit side.  The organising object is therefore the \emph{exit pair}, the joint law of $(T,J)$, or equivalently the two conditional laws $G_0=\Law(T\mid J=0)$ and $G_a=\Law(T\mid J=a)$ together with $\Pp_z(J=0)$.

Our first main result is a dictionary (\cref{thm:dictionary}).  Each classical stochastic order between $G_0$ and $G_a$ corresponds to one reset response: equality to universal neutrality, the usual order to a universal sign, the reversed hazard-rate order to monotonicity under deterministic deadlines, an order of generating functions to monotonicity under geometric resetting.  Structurally, the difference of the ruin probabilities under two clocks is the pairing on $\Lambda^2$ of the plane spanned by the clocks with the plane spanned by the exit pair (\cref{sec:plucker}).  The orders are then test regions for the Pl\"ucker coordinates $E_{nn'}=\nu_n\eta_{n'}-\nu_{n'}\eta_n$ of the exit pair, and the exact dual of hazard-rate ordered clocks turns out to be the reversed hazard-rate order (\cref{thm:exactdual}).

The rest of the paper decides which orders hold.  For even $a=2m$ a tail-swapping bijection identifies the relevant $2\times2$ determinants with weights of nonintersecting first-passage path pairs, which yields the likelihood-ratio order with an exact strictness classification (\cref{thm:discreteMain}).  For odd $a$ the two channels open on opposite parities, so $J$ is a function of $T$ and no site is universally neutral; nevertheless $\Law(T-1\mid J=a)$ dominates $\Law(T\mid J=0)$ in the likelihood-ratio order for $z<a/2$ (\cref{thm:odd}).  Consequently the responses that tolerate a one-step shift survive, together with monotonicity along clocks supported on one parity class (\cref{thm:oneparity}), while the pointwise ones fail; this explains the periodic-reset counterexample.  For general birth--death chains, without any symmetry, the direction of the exact-time response switches monotonically in the starting site and, when the second spectral mode is nondegenerate, for late times at the node of its eigenvector (\cref{thm:switch,prop:node}).  This spatial monotonicity is essentially the correct general statement: an exact rational example shows that a fixed site can reverse its temporal direction more than once (\cref{prop:singlefail}).  Finally, drifted Brownian motion shows the same architecture in the continuum, with an exact Poisson-reset response identity.

The companion paper \cite{PaperB} studies the complementary question of what complete reset neutrality, and the local response beyond it, reveal about the dynamics.

\section{Renewal reduction and exact-time exit information}

\subsection{Discrete renewal formula}
Consider the nearest-neighbour walk on $\{0,1,\ldots,a\}$ with
\[
\Pp(X_{n+1}=x+1\mid X_n=x)=p,\qquad
\Pp(X_{n+1}=x-1\mid X_n=x)=q,
\]
where $p,q>0$ and $p+q=1$.  Let
\[
T=\inf\{n\ge0:X_n\in\{0,a\}\},\qquad J=X_T\in\{0,a\}.
\]
For a fixed initial site $z\in\{1,\ldots,a-1\}$ define
\[
m_n=\Pp_z(T=n),\qquad
\nu_n=\Pp_z(T=n,J=0),\qquad
\eta_n=\Pp_z(T=n,J=a),
\]
and, whenever $m_n>0$,
\[
\varphi_n=\Pp_z(J=0\mid T=n)=\frac{\nu_n}{m_n}.
\]

Let $R$ be an independent reset clock.  We use the convention that absorption at step $n$ is counted whenever the walk is permitted to take that step.  Thus
\[
w_n=\Pp(R\ge n)
\]
is the survival weight of the reset clock.  After a reset the system restarts from the same $z$, independently of the past.  Throughout, a clock is called \emph{admissible} if $\sum_{n\ge1}m_nw_n>0$, so that an individual cycle has a positive probability of absorption.  All reset-response statements quantify only over admissible clocks, and deadline claims concern deadlines $N$ for which $\Pp_z(T\le N)>0$.  Without this restriction, a clock that resets before either boundary can be reached gives no eventual absorption and the renewal quotient is undefined.  The quotient is homogeneous in $w$, so proportional positive weight sequences represent the same reset response.

\begin{proposition}[Renewal weighted-average identity]\label[proposition]{prop:renewal}
For any admissible i.i.d. renewal resetting clock, the eventual lower-boundary absorption probability is
\begin{equation}
Q_z(w)=\frac{\sum_{n\ge1}\nu_n w_n}
{\sum_{n\ge1}m_n w_n}
=\frac{\sum_{n\ge1}\varphi_n m_n w_n}
{\sum_{n\ge1}m_n w_n}.
\label{eq:renewalaverage}
\end{equation}
\end{proposition}

\begin{proof}
Within one cycle let $A=\sum_n\nu_nw_n$ and $B=\sum_n\eta_nw_n$.  The remaining probability $1-A-B$ is reset without absorption.  Regeneration gives $Q=A+(1-A-B)Q$, hence $Q=A/(A+B)$, which is \eqref{eq:renewalaverage}.
\end{proof}

The interpretation is immediate: resetting does not create a new exit mechanism; it only reweights the absorption times.  This observation will be used repeatedly.  In particular, geometric survival weights are $s^{n-1}$ for a clock with $R\ge1$; the frequently used weights $s^n$ are proportional and yield exactly the same quotient.

\subsection{Bias cancellation}
Write $a=2m$ temporarily.  Let $N_-(n)$ and $N_+(n)$ be the unweighted numbers of admissible paths first exiting at time $n$ through $0$ and $2m$.  A path from $z$ to $0$ at time $n$ uses
\[
R=\frac{n-z}{2},\qquad L=\frac{n+z}{2},
\]
while a path from $z$ to $2m$ uses
\[
R=\frac{n+2m-z}{2},\qquad L=\frac{n-2m+z}{2}.
\]
Therefore
\begin{equation}
\frac{\eta_n}{\nu_n}
=\left(\frac pq\right)^m\frac{N_+(n)}{N_-(n)}.
\label{eq:biascancel}
\end{equation}
The dependence on $p/q$ is an $n$-independent factor.  Exact-time monotonicity is therefore a combinatorial first-passage property.

\subsection{Neutrality as independence}
The weighted-average identity gives a useful equivalence without assuming arbitrary independent weights.

\begin{theorem}[Universal neutrality equivalence]\label{thm:neutralind}
For a fixed reset site $z$, the following are equivalent.
\begin{enumerate}
\item $Q_z(w)$ is the same for every independent renewal reset law.
\item $Q_z(s)$ is the same for every geometric survival family $w_n=s^n$, $0<s<1$.
\item $\varphi_n$ is constant on the support of $T$.
\item $J$ and $T$ are independent under the reset-free law.
\end{enumerate}
\end{theorem}

\begin{proof}
Only $(2)\Rightarrow(3)$ needs comment.  If $Q_z(s)\equiv c$, then
\[
\sum_{n\ge1}(\nu_n-cm_n)s^n=0,\qquad 0<s<1.
\]
Analyticity of the power series gives $\nu_n=cm_n$ coefficientwise.  The equivalence of (3) and (4) is the definition of conditional independence.
\end{proof}

A bounded martingale formulation follows.  With $Y=\ind_{\{J=0\}}$ and $\mathcal G_n=\sigma(T\wedge n)$,
\[
M_n=\Ee[Y-\Ee Y\mid\mathcal G_n]
\]
is identically zero if and only if the site is universally reset-neutral.

\section{The exit pair and a dictionary of stochastic orders}\label{sec:dictionary}

By \cref{prop:renewal}, $Q_z(w)$ depends on the walk only through the joint law of $(T,J)$.  Write
\[
\pi_0=\Pp_z(J=0),\qquad
G_0=\Law(T\mid J=0),\qquad
G_a=\Law(T\mid J=a),
\]
and let $F_0,F_a$ be the corresponding distribution functions.  We call $(G_0,G_a;\pi_0)$ the \emph{exit pair}.  In this language
\[
Q_z(w)=\Pp(J=0\mid T\le R),
\]
with $R$ an independent copy of the reset clock.  A \emph{deadline} is the deterministic clock $R\equiv N$, i.e.\ $w_n=\ind_{\{n\le N\}}$, and we write $Q_z(N)$ for the corresponding ruin probability.

\begin{definition}\label[definition]{def:orders}
For laws $G,H$ on $\{1,2,\dots\}$ with mass functions $g,h$, distribution functions $F_G,F_H$ and generating functions $\phi_G,\phi_H$:
\begin{itemize}
\item $G\le_\st H$ if $F_G\ge F_H$;
\item $G\le_\rh H$ if $F_H/F_G$ is nondecreasing wherever defined, equivalently, if for all $N<M$,\[F_H(M)F_G(N)\ge F_H(N)F_G(M).\] This cross-product convention also covers zero values;
\item $G\le_\lr H$ if $g(n)h(n')\ge g(n')h(n)$ for all $n<n'$;
\item $G\le_\pgf H$ if $\phi_H(s)/\phi_G(s)$ is nondecreasing on $(0,1]$.
\end{itemize}
\end{definition}

Comparisons of nonnegative sequences under a nonincreasing ratio are understood in the cross-product sense $w_1(n)w_2(n')\ge w_1(n')w_2(n)$ for $n<n'$, even when some entries vanish.  The response quotients themselves still require admissibility.

It is classical that $\lr\Rightarrow\rh\Rightarrow\st$ \cite{ShakedShanthikumar2007}.  Also $\lr\Rightarrow\pgf$ (see the proof of \cref{thm:dictionary}).  The determinant form of $\le_\lr$ is essential, because it includes the times at which only one channel is open (\cref{rem:onechannel}).

\begin{theorem}[Dictionary]\label{thm:dictionary}
For every exit pair:
\begin{enumerate}[label=\textup{(\roman*)}]
\item \emph{Universal neutrality.} $Q_z(w)=\pi_0$ for every renewal clock if and only if $G_0=G_a$.
\item \emph{Universal sign.} $Q_z(w)\ge\pi_0$ for every renewal clock if and only if $G_0\le_\st G_a$.
\item \emph{Deadlines.} $N\mapsto Q_z(N)$ is nonincreasing if and only if $G_0\le_\rh G_a$.
\item \emph{Geometric resetting} (weights proportional to $s^n$). $s\mapsto Q_z(s)$ is nonincreasing on $(0,1]$ if and only if $G_0\le_\pgf G_a$.
\item \emph{Hazard-rate families.} If $G_0\le_\lr G_a$, then $Q_z(w_1)\ge Q_z(w_2)$ whenever $w_1/w_2$ is nonincreasing.
\end{enumerate}
Exchanging the roles of $0$ and $a$ reverses every statement.
\end{theorem}

\begin{proof}
For a deadline with $F_a(N)>0$,
\begin{equation}
\frac{Q_z(N)}{1-Q_z(N)}=\frac{\pi_0}{1-\pi_0}\,\frac{F_0(N)}{F_a(N)}.
\label{eq:deadlineodds}
\end{equation}
Item (iii) follows from \eqref{eq:deadlineodds} whenever both denominator factors are positive; at earlier admissible deadlines with $F_a(N)=0$, $Q_z(N)=1$, and the same assertion follows by cross multiplication.  For (ii), \eqref{eq:deadlineodds} gives $Q_z(N)\ge\pi_0$ if and only if $F_0(N)\ge F_a(N)$.  A general clock satisfies $w_n=\sum_N\Pp(R=N)\ind_{\{n\le N\}}$, so numerator and denominator of \eqref{eq:renewalaverage} are mixtures of those of the deadlines, and $Q_z(w)$ is a mediant of the values $Q_z(N)$.  Item (i) is (ii) applied in both directions (\cref{thm:neutralind} shows that the geometric family alone already suffices).  Item (iv) is the analogue of \eqref{eq:deadlineodds} with $w_n=s^n$.  For (v), $Q_z(w_1)\ge Q_z(w_2)$ is equivalent to $\sum_{n,n'}w_1(n)w_2(n')\big(\nu_n\eta_{n'}-\eta_n\nu_{n'}\big)\ge0$.  Pairing $(n,n')$ with $(n',n)$, the sum equals
\[
\sum_{n<n'}\big(w_1(n)w_2(n')-w_1(n')w_2(n)\big)\big(\nu_n\eta_{n'}-\nu_{n'}\eta_n\big),
\]
where the first factor is $\ge0$ because $w_1/w_2$ is nonincreasing, and the second is $\ge0$ by $G_0\le_\lr G_a$.  The same computation with $w_1=s^n$, $w_2=t^n$, $s<t$, gives $\lr\Rightarrow\pgf$.
\end{proof}

\begin{remark}[Shifts]\label[remark]{rem:shift}
If $\Law(T-d\mid J=a)\ge_\lr G_0$ for some $d\ge0$, then $G_0\le_\st G_a$ and $G_0\le_\pgf G_a$ (the shift multiplies $\phi_{G_a}$ by $s^d$), but in general neither (iii) nor (v) follows.  This is the mechanism behind the odd case (\cref{sec:odd}).
\end{remark}

\begin{remark}[Times with one open channel]\label[remark]{rem:onechannel}
The times at which only one side can be reached carry weight in \eqref{eq:renewalaverage} and must be kept.  Discarding them produces false neutrality: for $a=4$, $z=1$ the exact-time odds restricted to times with both channels open are constant, yet $Q_z$ depends on the clock.  All order statements are therefore written as determinant inequalities, which include these times.
\end{remark}

\subsection{Pl\"ucker formulation and exact duality}\label{sec:plucker}

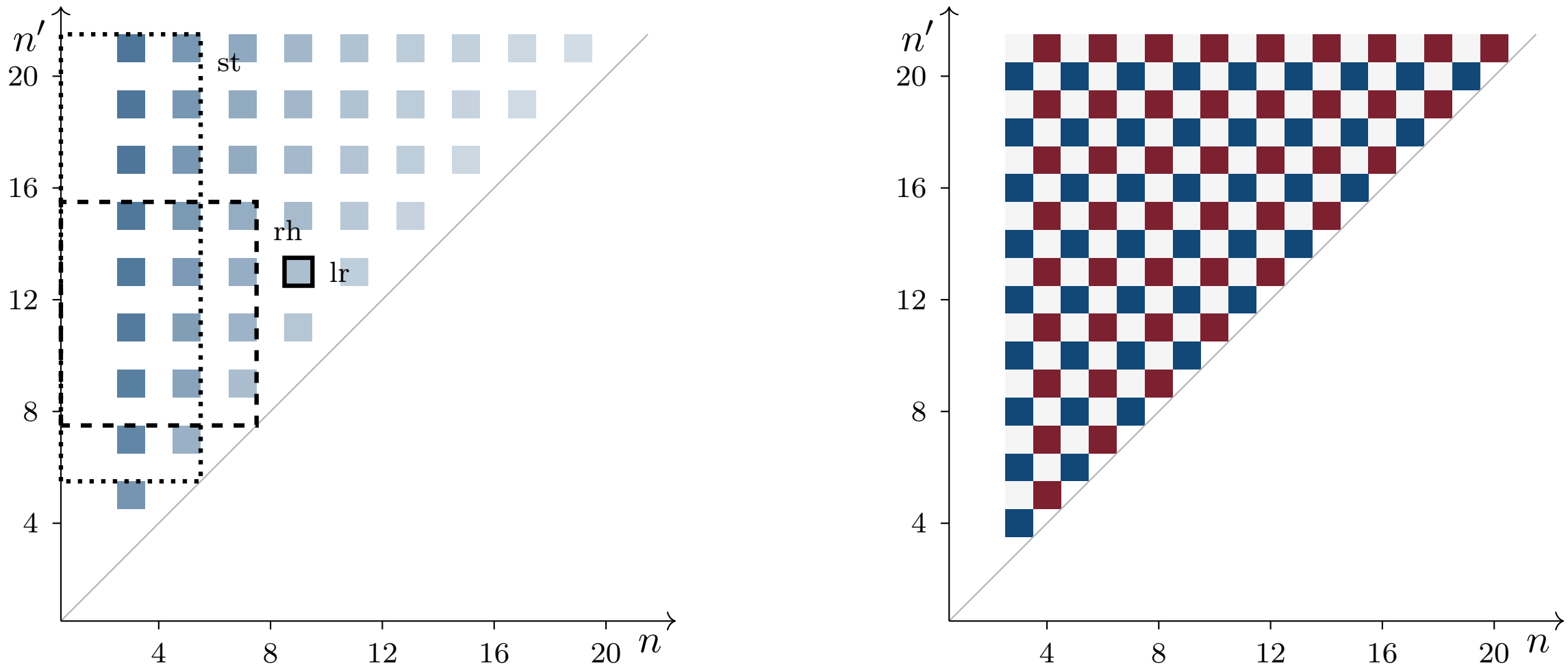
\begin{figure}[t]
\centering
\begin{subfigure}[b]{0.48\textwidth}\centering
\begin{tikzpicture}[x=0.27cm,y=0.27cm]
\fill[deepblue!58!white] (2.5,4.5) rectangle (3.5,5.5);
\fill[deepblue!66!white] (2.5,6.5) rectangle (3.5,7.5);
\fill[deepblue!70!white] (2.5,8.5) rectangle (3.5,9.5);
\fill[deepblue!72!white] (2.5,10.5) rectangle (3.5,11.5);
\fill[deepblue!73!white] (2.5,12.5) rectangle (3.5,13.5);
\fill[deepblue!74!white] (2.5,14.5) rectangle (3.5,15.5);
\fill[deepblue!75!white] (2.5,16.5) rectangle (3.5,17.5);
\fill[deepblue!75!white] (2.5,18.5) rectangle (3.5,19.5);
\fill[deepblue!75!white] (2.5,20.5) rectangle (3.5,21.5);
\fill[deepblue!43!white] (4.5,6.5) rectangle (5.5,7.5);
\fill[deepblue!50!white] (4.5,8.5) rectangle (5.5,9.5);
\fill[deepblue!53!white] (4.5,10.5) rectangle (5.5,11.5);
\fill[deepblue!55!white] (4.5,12.5) rectangle (5.5,13.5);
\fill[deepblue!56!white] (4.5,14.5) rectangle (5.5,15.5);
\fill[deepblue!57!white] (4.5,16.5) rectangle (5.5,17.5);
\fill[deepblue!57!white] (4.5,18.5) rectangle (5.5,19.5);
\fill[deepblue!57!white] (4.5,20.5) rectangle (5.5,21.5);
\fill[deepblue!36!white] (6.5,8.5) rectangle (7.5,9.5);
\fill[deepblue!41!white] (6.5,10.5) rectangle (7.5,11.5);
\fill[deepblue!44!white] (6.5,12.5) rectangle (7.5,13.5);
\fill[deepblue!45!white] (6.5,14.5) rectangle (7.5,15.5);
\fill[deepblue!46!white] (6.5,16.5) rectangle (7.5,17.5);
\fill[deepblue!46!white] (6.5,18.5) rectangle (7.5,19.5);
\fill[deepblue!47!white] (6.5,20.5) rectangle (7.5,21.5);
\fill[deepblue!31!white] (8.5,10.5) rectangle (9.5,11.5);
\fill[deepblue!35!white] (8.5,12.5) rectangle (9.5,13.5);
\fill[deepblue!37!white] (8.5,14.5) rectangle (9.5,15.5);
\fill[deepblue!38!white] (8.5,16.5) rectangle (9.5,17.5);
\fill[deepblue!39!white] (8.5,18.5) rectangle (9.5,19.5);
\fill[deepblue!39!white] (8.5,20.5) rectangle (9.5,21.5);
\fill[deepblue!27!white] (10.5,12.5) rectangle (11.5,13.5);
\fill[deepblue!30!white] (10.5,14.5) rectangle (11.5,15.5);
\fill[deepblue!32!white] (10.5,16.5) rectangle (11.5,17.5);
\fill[deepblue!33!white] (10.5,18.5) rectangle (11.5,19.5);
\fill[deepblue!33!white] (10.5,20.5) rectangle (11.5,21.5);
\fill[deepblue!24!white] (12.5,14.5) rectangle (13.5,15.5);
\fill[deepblue!27!white] (12.5,16.5) rectangle (13.5,17.5);
\fill[deepblue!28!white] (12.5,18.5) rectangle (13.5,19.5);
\fill[deepblue!28!white] (12.5,20.5) rectangle (13.5,21.5);
\fill[deepblue!22!white] (14.5,16.5) rectangle (15.5,17.5);
\fill[deepblue!24!white] (14.5,18.5) rectangle (15.5,19.5);
\fill[deepblue!25!white] (14.5,20.5) rectangle (15.5,21.5);
\fill[deepblue!20!white] (16.5,18.5) rectangle (17.5,19.5);
\fill[deepblue!22!white] (16.5,20.5) rectangle (17.5,21.5);
\fill[deepblue!19!white] (18.5,20.5) rectangle (19.5,21.5);
\draw[gray!60] (0.5,0.5) -- (21.5,21.5);
\draw[->] (0.5,0.5) -- (22.5,0.5) node[below left] {$n$};
\draw[->] (0.5,0.5) -- (0.5,22.5) node[below left] {$n'$};
\draw (4,0.5) -- (4,0.2) node[below] {\scriptsize 4};
\draw (0.5,4) -- (0.2,4) node[left] {\scriptsize 4};
\draw (8,0.5) -- (8,0.2) node[below] {\scriptsize 8};
\draw (0.5,8) -- (0.2,8) node[left] {\scriptsize 8};
\draw (12,0.5) -- (12,0.2) node[below] {\scriptsize 12};
\draw (0.5,12) -- (0.2,12) node[left] {\scriptsize 12};
\draw (16,0.5) -- (16,0.2) node[below] {\scriptsize 16};
\draw (0.5,16) -- (0.2,16) node[left] {\scriptsize 16};
\draw (20,0.5) -- (20,0.2) node[below] {\scriptsize 20};
\draw (0.5,20) -- (0.2,20) node[left] {\scriptsize 20};
\draw[very thick,dotted] (0.5,5.5) rectangle (5.5,21.5);
\node[font=\scriptsize,anchor=west] at (5.6,20.5) {st};
\draw[very thick,dashed] (0.5,7.5) rectangle (7.5,15.5);
\node[font=\scriptsize,anchor=west] at (7.6,14.5) {rh};
\draw[very thick] (8.5,12.5) rectangle (9.5,13.5);
\node[font=\scriptsize,anchor=west] at (9.6,13) {lr};
\end{tikzpicture}
\caption{$a=8$, $z=3$: $E\ge0$ everywhere ($\lr$).}\end{subfigure}\hfill
\begin{subfigure}[b]{0.48\textwidth}\centering
\begin{tikzpicture}[x=0.27cm,y=0.27cm]
\fill[deepblue!100!white] (2.5,3.5) rectangle (3.5,4.5);
\fill[gray!8] (2.5,4.5) rectangle (3.5,5.5);
\fill[deepblue!100!white] (2.5,5.5) rectangle (3.5,6.5);
\fill[gray!8] (2.5,6.5) rectangle (3.5,7.5);
\fill[deepblue!100!white] (2.5,7.5) rectangle (3.5,8.5);
\fill[gray!8] (2.5,8.5) rectangle (3.5,9.5);
\fill[deepblue!100!white] (2.5,9.5) rectangle (3.5,10.5);
\fill[gray!8] (2.5,10.5) rectangle (3.5,11.5);
\fill[deepblue!100!white] (2.5,11.5) rectangle (3.5,12.5);
\fill[gray!8] (2.5,12.5) rectangle (3.5,13.5);
\fill[deepblue!100!white] (2.5,13.5) rectangle (3.5,14.5);
\fill[gray!8] (2.5,14.5) rectangle (3.5,15.5);
\fill[deepblue!100!white] (2.5,15.5) rectangle (3.5,16.5);
\fill[gray!8] (2.5,16.5) rectangle (3.5,17.5);
\fill[deepblue!100!white] (2.5,17.5) rectangle (3.5,18.5);
\fill[gray!8] (2.5,18.5) rectangle (3.5,19.5);
\fill[deepblue!100!white] (2.5,19.5) rectangle (3.5,20.5);
\fill[gray!8] (2.5,20.5) rectangle (3.5,21.5);
\fill[burgundy!100!white] (3.5,4.5) rectangle (4.5,5.5);
\fill[gray!8] (3.5,5.5) rectangle (4.5,6.5);
\fill[burgundy!100!white] (3.5,6.5) rectangle (4.5,7.5);
\fill[gray!8] (3.5,7.5) rectangle (4.5,8.5);
\fill[burgundy!100!white] (3.5,8.5) rectangle (4.5,9.5);
\fill[gray!8] (3.5,9.5) rectangle (4.5,10.5);
\fill[burgundy!100!white] (3.5,10.5) rectangle (4.5,11.5);
\fill[gray!8] (3.5,11.5) rectangle (4.5,12.5);
\fill[burgundy!100!white] (3.5,12.5) rectangle (4.5,13.5);
\fill[gray!8] (3.5,13.5) rectangle (4.5,14.5);
\fill[burgundy!100!white] (3.5,14.5) rectangle (4.5,15.5);
\fill[gray!8] (3.5,15.5) rectangle (4.5,16.5);
\fill[burgundy!100!white] (3.5,16.5) rectangle (4.5,17.5);
\fill[gray!8] (3.5,17.5) rectangle (4.5,18.5);
\fill[burgundy!100!white] (3.5,18.5) rectangle (4.5,19.5);
\fill[gray!8] (3.5,19.5) rectangle (4.5,20.5);
\fill[burgundy!100!white] (3.5,20.5) rectangle (4.5,21.5);
\fill[deepblue!100!white] (4.5,5.5) rectangle (5.5,6.5);
\fill[gray!8] (4.5,6.5) rectangle (5.5,7.5);
\fill[deepblue!100!white] (4.5,7.5) rectangle (5.5,8.5);
\fill[gray!8] (4.5,8.5) rectangle (5.5,9.5);
\fill[deepblue!100!white] (4.5,9.5) rectangle (5.5,10.5);
\fill[gray!8] (4.5,10.5) rectangle (5.5,11.5);
\fill[deepblue!100!white] (4.5,11.5) rectangle (5.5,12.5);
\fill[gray!8] (4.5,12.5) rectangle (5.5,13.5);
\fill[deepblue!100!white] (4.5,13.5) rectangle (5.5,14.5);
\fill[gray!8] (4.5,14.5) rectangle (5.5,15.5);
\fill[deepblue!100!white] (4.5,15.5) rectangle (5.5,16.5);
\fill[gray!8] (4.5,16.5) rectangle (5.5,17.5);
\fill[deepblue!100!white] (4.5,17.5) rectangle (5.5,18.5);
\fill[gray!8] (4.5,18.5) rectangle (5.5,19.5);
\fill[deepblue!100!white] (4.5,19.5) rectangle (5.5,20.5);
\fill[gray!8] (4.5,20.5) rectangle (5.5,21.5);
\fill[burgundy!100!white] (5.5,6.5) rectangle (6.5,7.5);
\fill[gray!8] (5.5,7.5) rectangle (6.5,8.5);
\fill[burgundy!100!white] (5.5,8.5) rectangle (6.5,9.5);
\fill[gray!8] (5.5,9.5) rectangle (6.5,10.5);
\fill[burgundy!100!white] (5.5,10.5) rectangle (6.5,11.5);
\fill[gray!8] (5.5,11.5) rectangle (6.5,12.5);
\fill[burgundy!100!white] (5.5,12.5) rectangle (6.5,13.5);
\fill[gray!8] (5.5,13.5) rectangle (6.5,14.5);
\fill[burgundy!100!white] (5.5,14.5) rectangle (6.5,15.5);
\fill[gray!8] (5.5,15.5) rectangle (6.5,16.5);
\fill[burgundy!100!white] (5.5,16.5) rectangle (6.5,17.5);
\fill[gray!8] (5.5,17.5) rectangle (6.5,18.5);
\fill[burgundy!100!white] (5.5,18.5) rectangle (6.5,19.5);
\fill[gray!8] (5.5,19.5) rectangle (6.5,20.5);
\fill[burgundy!100!white] (5.5,20.5) rectangle (6.5,21.5);
\fill[deepblue!100!white] (6.5,7.5) rectangle (7.5,8.5);
\fill[gray!8] (6.5,8.5) rectangle (7.5,9.5);
\fill[deepblue!100!white] (6.5,9.5) rectangle (7.5,10.5);
\fill[gray!8] (6.5,10.5) rectangle (7.5,11.5);
\fill[deepblue!100!white] (6.5,11.5) rectangle (7.5,12.5);
\fill[gray!8] (6.5,12.5) rectangle (7.5,13.5);
\fill[deepblue!100!white] (6.5,13.5) rectangle (7.5,14.5);
\fill[gray!8] (6.5,14.5) rectangle (7.5,15.5);
\fill[deepblue!100!white] (6.5,15.5) rectangle (7.5,16.5);
\fill[gray!8] (6.5,16.5) rectangle (7.5,17.5);
\fill[deepblue!100!white] (6.5,17.5) rectangle (7.5,18.5);
\fill[gray!8] (6.5,18.5) rectangle (7.5,19.5);
\fill[deepblue!100!white] (6.5,19.5) rectangle (7.5,20.5);
\fill[gray!8] (6.5,20.5) rectangle (7.5,21.5);
\fill[burgundy!100!white] (7.5,8.5) rectangle (8.5,9.5);
\fill[gray!8] (7.5,9.5) rectangle (8.5,10.5);
\fill[burgundy!100!white] (7.5,10.5) rectangle (8.5,11.5);
\fill[gray!8] (7.5,11.5) rectangle (8.5,12.5);
\fill[burgundy!100!white] (7.5,12.5) rectangle (8.5,13.5);
\fill[gray!8] (7.5,13.5) rectangle (8.5,14.5);
\fill[burgundy!100!white] (7.5,14.5) rectangle (8.5,15.5);
\fill[gray!8] (7.5,15.5) rectangle (8.5,16.5);
\fill[burgundy!100!white] (7.5,16.5) rectangle (8.5,17.5);
\fill[gray!8] (7.5,17.5) rectangle (8.5,18.5);
\fill[burgundy!100!white] (7.5,18.5) rectangle (8.5,19.5);
\fill[gray!8] (7.5,19.5) rectangle (8.5,20.5);
\fill[burgundy!100!white] (7.5,20.5) rectangle (8.5,21.5);
\fill[deepblue!100!white] (8.5,9.5) rectangle (9.5,10.5);
\fill[gray!8] (8.5,10.5) rectangle (9.5,11.5);
\fill[deepblue!100!white] (8.5,11.5) rectangle (9.5,12.5);
\fill[gray!8] (8.5,12.5) rectangle (9.5,13.5);
\fill[deepblue!100!white] (8.5,13.5) rectangle (9.5,14.5);
\fill[gray!8] (8.5,14.5) rectangle (9.5,15.5);
\fill[deepblue!100!white] (8.5,15.5) rectangle (9.5,16.5);
\fill[gray!8] (8.5,16.5) rectangle (9.5,17.5);
\fill[deepblue!100!white] (8.5,17.5) rectangle (9.5,18.5);
\fill[gray!8] (8.5,18.5) rectangle (9.5,19.5);
\fill[deepblue!100!white] (8.5,19.5) rectangle (9.5,20.5);
\fill[gray!8] (8.5,20.5) rectangle (9.5,21.5);
\fill[burgundy!100!white] (9.5,10.5) rectangle (10.5,11.5);
\fill[gray!8] (9.5,11.5) rectangle (10.5,12.5);
\fill[burgundy!100!white] (9.5,12.5) rectangle (10.5,13.5);
\fill[gray!8] (9.5,13.5) rectangle (10.5,14.5);
\fill[burgundy!100!white] (9.5,14.5) rectangle (10.5,15.5);
\fill[gray!8] (9.5,15.5) rectangle (10.5,16.5);
\fill[burgundy!100!white] (9.5,16.5) rectangle (10.5,17.5);
\fill[gray!8] (9.5,17.5) rectangle (10.5,18.5);
\fill[burgundy!100!white] (9.5,18.5) rectangle (10.5,19.5);
\fill[gray!8] (9.5,19.5) rectangle (10.5,20.5);
\fill[burgundy!100!white] (9.5,20.5) rectangle (10.5,21.5);
\fill[deepblue!100!white] (10.5,11.5) rectangle (11.5,12.5);
\fill[gray!8] (10.5,12.5) rectangle (11.5,13.5);
\fill[deepblue!100!white] (10.5,13.5) rectangle (11.5,14.5);
\fill[gray!8] (10.5,14.5) rectangle (11.5,15.5);
\fill[deepblue!100!white] (10.5,15.5) rectangle (11.5,16.5);
\fill[gray!8] (10.5,16.5) rectangle (11.5,17.5);
\fill[deepblue!100!white] (10.5,17.5) rectangle (11.5,18.5);
\fill[gray!8] (10.5,18.5) rectangle (11.5,19.5);
\fill[deepblue!100!white] (10.5,19.5) rectangle (11.5,20.5);
\fill[gray!8] (10.5,20.5) rectangle (11.5,21.5);
\fill[burgundy!100!white] (11.5,12.5) rectangle (12.5,13.5);
\fill[gray!8] (11.5,13.5) rectangle (12.5,14.5);
\fill[burgundy!100!white] (11.5,14.5) rectangle (12.5,15.5);
\fill[gray!8] (11.5,15.5) rectangle (12.5,16.5);
\fill[burgundy!100!white] (11.5,16.5) rectangle (12.5,17.5);
\fill[gray!8] (11.5,17.5) rectangle (12.5,18.5);
\fill[burgundy!100!white] (11.5,18.5) rectangle (12.5,19.5);
\fill[gray!8] (11.5,19.5) rectangle (12.5,20.5);
\fill[burgundy!100!white] (11.5,20.5) rectangle (12.5,21.5);
\fill[deepblue!100!white] (12.5,13.5) rectangle (13.5,14.5);
\fill[gray!8] (12.5,14.5) rectangle (13.5,15.5);
\fill[deepblue!100!white] (12.5,15.5) rectangle (13.5,16.5);
\fill[gray!8] (12.5,16.5) rectangle (13.5,17.5);
\fill[deepblue!100!white] (12.5,17.5) rectangle (13.5,18.5);
\fill[gray!8] (12.5,18.5) rectangle (13.5,19.5);
\fill[deepblue!100!white] (12.5,19.5) rectangle (13.5,20.5);
\fill[gray!8] (12.5,20.5) rectangle (13.5,21.5);
\fill[burgundy!100!white] (13.5,14.5) rectangle (14.5,15.5);
\fill[gray!8] (13.5,15.5) rectangle (14.5,16.5);
\fill[burgundy!100!white] (13.5,16.5) rectangle (14.5,17.5);
\fill[gray!8] (13.5,17.5) rectangle (14.5,18.5);
\fill[burgundy!100!white] (13.5,18.5) rectangle (14.5,19.5);
\fill[gray!8] (13.5,19.5) rectangle (14.5,20.5);
\fill[burgundy!100!white] (13.5,20.5) rectangle (14.5,21.5);
\fill[deepblue!100!white] (14.5,15.5) rectangle (15.5,16.5);
\fill[gray!8] (14.5,16.5) rectangle (15.5,17.5);
\fill[deepblue!100!white] (14.5,17.5) rectangle (15.5,18.5);
\fill[gray!8] (14.5,18.5) rectangle (15.5,19.5);
\fill[deepblue!100!white] (14.5,19.5) rectangle (15.5,20.5);
\fill[gray!8] (14.5,20.5) rectangle (15.5,21.5);
\fill[burgundy!100!white] (15.5,16.5) rectangle (16.5,17.5);
\fill[gray!8] (15.5,17.5) rectangle (16.5,18.5);
\fill[burgundy!100!white] (15.5,18.5) rectangle (16.5,19.5);
\fill[gray!8] (15.5,19.5) rectangle (16.5,20.5);
\fill[burgundy!100!white] (15.5,20.5) rectangle (16.5,21.5);
\fill[deepblue!100!white] (16.5,17.5) rectangle (17.5,18.5);
\fill[gray!8] (16.5,18.5) rectangle (17.5,19.5);
\fill[deepblue!100!white] (16.5,19.5) rectangle (17.5,20.5);
\fill[gray!8] (16.5,20.5) rectangle (17.5,21.5);
\fill[burgundy!100!white] (17.5,18.5) rectangle (18.5,19.5);
\fill[gray!8] (17.5,19.5) rectangle (18.5,20.5);
\fill[burgundy!100!white] (17.5,20.5) rectangle (18.5,21.5);
\fill[deepblue!100!white] (18.5,19.5) rectangle (19.5,20.5);
\fill[gray!8] (18.5,20.5) rectangle (19.5,21.5);
\fill[burgundy!100!white] (19.5,20.5) rectangle (20.5,21.5);
\draw[gray!60] (0.5,0.5) -- (21.5,21.5);
\draw[->] (0.5,0.5) -- (22.5,0.5) node[below left] {$n$};
\draw[->] (0.5,0.5) -- (0.5,22.5) node[below left] {$n'$};
\draw (4,0.5) -- (4,0.2) node[below] {\scriptsize 4};
\draw (0.5,4) -- (0.2,4) node[left] {\scriptsize 4};
\draw (8,0.5) -- (8,0.2) node[below] {\scriptsize 8};
\draw (0.5,8) -- (0.2,8) node[left] {\scriptsize 8};
\draw (12,0.5) -- (12,0.2) node[below] {\scriptsize 12};
\draw (0.5,12) -- (0.2,12) node[left] {\scriptsize 12};
\draw (16,0.5) -- (16,0.2) node[below] {\scriptsize 16};
\draw (0.5,16) -- (0.2,16) node[left] {\scriptsize 16};
\draw (20,0.5) -- (20,0.2) node[below] {\scriptsize 20};
\draw (0.5,20) -- (0.2,20) node[left] {\scriptsize 20};
\end{tikzpicture}
\caption{$a=7$, $z=3$: parity checkerboard.}\end{subfigure}
\caption{Pl\"ucker coordinates of the exit pair for the symmetric walk, on the triangle $n<n'$.  Each cell shows the normalised coordinate $E_{nn'}/(m_nm_{n'})=\varphi_n-\varphi_{n'}$, in blue when positive and red when negative, with intensity proportional to the square root of its absolute value; light grey cells are zero, and blank cells are times at which no absorption is possible.  In (a) the test regions of the orders are marked: a single cell ($\lr$), an anchored rectangle $[1,N]\times(N,M]$ ($\rh$) and the displayed truncation of an infinite anchored rectangle $[1,N]\times(N,\infty)$ ($\st$).  In (b) the two channels open on opposite parities, so the pointwise orders fail, while the pair shifted by one step is again likelihood-ratio ordered (\cref{thm:odd}).  The displayed panels are generated by \texttt{plucker\_figure.c} using floating-point forward recursion; the panel signs were also independently checked in exact rational arithmetic.}
\label{fig:plucker}
\end{figure}

Let $m_n=\nu_n+\eta_n$ and write $\langle w,x\rangle=\sum_nw_nx_n$.  For two clocks,
\begin{equation}
Q_z(w_1)-Q_z(w_2)=\frac{\langle w_1\wedge w_2,\ \nu\wedge\eta\rangle}{\langle w_1,m\rangle\,\langle w_2,m\rangle},
\qquad
\langle a\wedge b,c\wedge d\rangle=\langle a,c\rangle\langle b,d\rangle-\langle a,d\rangle\langle b,c\rangle .
\notag
\end{equation}
Comparing two reset mechanisms is therefore the canonical pairing, on $\Lambda^2$ of the sequence space, between the plane spanned by the two clocks and the plane spanned by the exit pair.  By Cauchy--Binet the pairing is $\sum_{n<n'}(w_1\wedge w_2)_{nn'}E_{nn'}$, where
\[
E_{nn'}=\nu_n\eta_{n'}-\nu_{n'}\eta_n,\qquad n<n',
\]
are the Pl\"ucker coordinates of the exit pair.  In these coordinates:
\begin{itemize}
\item universal neutrality is $E=0$ (the plane degenerates to a line, $\nu\parallel\eta$);
\item $G_0\le_\lr G_a$ is $E\ge0$ pointwise on the triangle $n<n'$;
\item $G_0\le_\rh G_a$ is $\sum E\ge0$ over every rectangle $[1,N]\times(N,M]$, since for deadlines $d_N=\ind_{[1,N]}$ one has $d_N\wedge d_M=\sum_{n\le N<n'\le M}e_n\wedge e_{n'}$;
\item $G_0\le_\st G_a$ is $\sum E\ge0$ over every infinite rectangle $[1,N]\times(N,\infty)$;
\item $G_0\le_\pgf G_a$ is $\sum E_{nn'}(s^nt^{n'}-s^{n'}t^n)\ge0$ for $s<t$.
\end{itemize}
The hierarchy of stochastic orders is thus a hierarchy of test regions for one density on the triangle.  This is the dual-cone description familiar from integral stochastic orders \cite{MullerStoyan2002}, in which an order is generated by a cone of test functions; here the test functions are the antisymmetric products $w_1\wedge w_2$ of reset clocks.  The covariance identity of \cref{thm:covresp} is the infinitesimal version: $\dd Q_z=-\langle w\wedge\partial_\theta w,\ \nu\wedge\eta\rangle\,\dd\theta/\langle w,m\rangle^2$.  The path-switching lemma gives the Pl\"ucker coordinates a combinatorial meaning: after reflection, $E_{nn'}$ is the weight of nonintersecting path pairs, in the spirit of Lindstr\"om--Gessel--Viennot \cite{Lindstrom1973,GesselViennot1985}.  \Cref{fig:plucker} shows the density and the test regions.

\begin{theorem}[Exact duality]\label{thm:exactdual}
For every exit pair the following are equivalent:
\begin{enumerate}[label=\textup{(\alph*)}]
\item $G_0\le_\rh G_a$;
\item $Q_z(w_1)\ge Q_z(w_2)$ for every pair of clocks with $w_1/w_2$ nonincreasing (hazard-rate ordered, clock~1 more intense);
\item $Q_z(R_1)\ge Q_z(R_2)$ for every pair of clocks with $R_1\le_\lr R_2$.
\end{enumerate}
Moreover, $G_0\le_\lr G_a$ holds if and only if $\langle w_1\wedge w_2,\nu\wedge\eta\rangle\ge0$ (equivalently $Q_z(w_1)\ge Q_z(w_2)$ whenever both quotients are defined) for all nonnegative weight sequences, not necessarily survival functions, with $w_1/w_2$ nonincreasing.
\end{theorem}

\begin{proof}
(b)$\Rightarrow$(a) and (c)$\Rightarrow$(a): deadlines $N<M$ are both hazard-rate and likelihood-ratio ordered, and by \cref{thm:dictionary}(iii) monotone deadline response is (a).

(a)$\Rightarrow$(c): with $a_N=\Pp(R_1=N)$ and $b_N=\Pp(R_2=N)$, bilinearity gives $w_1\wedge w_2=\sum_{N<M}(a_Nb_M-a_Mb_N)\,d_N\wedge d_M$, with nonnegative coefficients when $R_1\le_\lr R_2$.  Moreover $\langle d_N\wedge d_M,\nu\wedge\eta\rangle=B_NA_M-B_MA_N$, with $B_N=\sum_{n\le N}\nu_n$ and $A_N=\sum_{n\le N}\eta_n$, which is $\ge0$ by (a).

(a)$\Rightarrow$(b): clocks satisfy $R\ge1$, so $w_1(1)=w_2(1)=1$.  Write $w_1=hw_2$ with $h$ nonincreasing and $h(1)=1$ (on any tail where $w_2$ vanishes, $w_1$ vanishes as well and $h$ may be set to $0$), so $h$ is itself a survival function.  Then $Q_z(w_1)\ge Q_z(w_2)$ is the universal-sign inequality, with clock $h$, for the tilted pair $(w_2\nu,w_2\eta)$; by \cref{thm:dictionary}(ii) it suffices that the tilted pair be $\st$-ordered.  The reversed hazard-rate order is preserved by nonincreasing tilts.  Preservation of stochastic orders under a common weighting is a classical theme \cite{ShakedShanthikumar2007,BartoszewiczSkolimowska2006}; we include the short argument in the form needed here.  To see this, write $w_2(n)=\sum_{k\ge n}c_k$ with $c_k\ge0$ (including a possible atom at $k=\infty$) and put $\varrho(N)=F_a(N)/F_0(N)$ wherever defined.  The unnormalised tilted cumulative masses at $N$ are $H_i(N)=\sum_k c_k F_i(k\wedge N)$ for $i=0,a$.  Passing from $N$ to $N+1$ increases them by $C_N[F_i(N+1)-F_i(N)]$, where $C_N=\sum_{k>N}c_k\ge0$.  Since $\varrho$ is nondecreasing, the ratio of these increments (when the denominator is positive) is at least $\varrho(N)$; in turn $H_a(N)/H_0(N)$ is a weighted mean of $\varrho(k\wedge N)$ and is at most $\varrho(N)$.  The mediant inequality therefore gives $H_a(N+1)/H_0(N+1)\ge H_a(N)/H_0(N)$.  Zero denominators follow directly by cross multiplication.  Hence the tilted pair is $\rh$-, and therefore $\st$-ordered.

For the last statement, sufficiency is the determinant computation in the proof of \cref{thm:dictionary}(v); for necessity take $w_1=e_n$, $w_2=e_{n'}$ with $n<n'$, so that $w_1\wedge w_2=e_n\wedge e_{n'}$ and the pairing is $E_{nn'}$; the pairing form avoids any condition $m_n>0$.
\end{proof}

\begin{remark}
Thus reversed hazard-rate order, which is equivalent to monotone response under deadlines, is exactly what renewal resetting along hazard-rate or likelihood-ratio ordered clocks requires.  The full likelihood-ratio order of \cref{cor:evenlr} corresponds to a stronger, gated response: arbitrary nonnegative detection weights rather than survival functions.
\end{remark}

\begin{proposition}[Usual stochastic order of clocks is insufficient]\label[proposition]{prop:clockstfail}
Even for the symmetric nearest-neighbour walk, $R_1\le_\st R_2$ does not imply $Q_z(R_1)\ge Q_z(R_2)$.  For the walk on $\{0,\ldots,6\}$ started at $z=1$, let
\[
\Pp(R_1=1)=\Pp(R_1=5)=\frac12,
\qquad
\Pp(R_2=3)=\Pp(R_2=5)=\frac12.
\]
Then $R_1\le_\st R_2$, but
\[
Q_1(R_1)=\frac{38}{39}<\frac{42}{43}=Q_1(R_2).
\]
\end{proposition}

\begin{proof}
Up to time $5$ the lower-channel masses are $\nu_1=1/2$, $\nu_3=1/8$, $\nu_5=1/16$, while the only upper-channel mass is $\eta_5=1/32$.  For $R_1$, the survival weights are $1$ at time $1$ and $1/2$ at times $3,5$, giving weighted lower/upper masses $19/32$ and $1/64$.  For $R_2$, the weights are $1$ at times $1,3$ and $1/2$ at time $5$, giving $21/32$ and $1/64$.  The displayed ratios follow.  The stochastic order is immediate because the atom of size $1/2$ is shifted from time $1$ to time $3$ while the atom at time $5$ is unchanged.
\end{proof}

\section{Even intervals: corrected monotonicity and a path-switching proof}

\subsection{Parity and the corrected statement}
Suppose $a=2m$.  Both boundaries can be reached only on the parity class $n\equiv z\pmod2$.  If $z<m$, the lower boundary is reachable from time $z$, whereas the upper boundary cannot be reached before time $2m-z$.  Hence
\begin{equation}
\varphi_n=1,
\qquad z\le n<2m-z,
\quad n\equiv z\pmod2.
\notag
\end{equation}
The strict statement on the entire support is therefore impossible.  The correct theorem separates the one-channel plateau from the two-channel regime.

\begin{theorem}[Corrected even-interval exit-time monotonicity]\label{thm:discreteMain}
Let $a=2m$ and $0<z<2m$.
\begin{enumerate}
\item If $z<m$, then $\varphi_n$ is nonincreasing on its admissible parity.  It equals $1$ before the upper channel opens.  If $a\ge6$, it is strictly decreasing at every comparison after both channels are open.
\item If $z>m$, then $\varphi_n$ is nondecreasing, with the reflected strictness statement.
\item If $z=m$, then
\begin{equation}
\varphi_n\equiv\varphi_*=\frac{q^m}{p^m+q^m}.
\notag
\end{equation}
\item The only off-midpoint failure of post-opening strictness is $a=4,z=1$ and its reflection $z=3$.  For $a=4,z=1$,
\[
\varphi_1=1,
\qquad
\varphi_n=\frac{q^2}{p^2+q^2},\qquad n=3,5,7,\ldots.
\]
\item For every $z$,
\begin{equation}
\varphi_n\longrightarrow\varphi_*
\qquad(n\to\infty\text{ on the admissible parity}).
\label{eq:philimit}
\end{equation}
\end{enumerate}
\end{theorem}

The key step is a first-passage $TP_2$ inequality.  Classical total-positivity results for absorption arrays go back to Karlin and Karlin--McGregor \cite{Karlin1964,KarlinMcGregor1960}.  Here the order-two statement needed for \cref{thm:discreteMain} admits a short direct proof.

\subsection{First-passage tail switching}
For states $i\in\{1,\ldots,2m-1\}$, let $F_n(i)$ be the total weight of paths which start at $i$ and first hit $0$ at time $n$ before hitting $2m$.

\begin{lemma}[Path-switching determinant]\label[lemma]{lem:pathswitch}
Let $i<j$ have the same parity and let $r<s$ be admissible first-passage times.  Then
\begin{equation}
F_r(i)F_s(j)-F_r(j)F_s(i)\ge0.
\label{eq:TP2}
\end{equation}
Moreover the left side equals the total pair-weight of nonintersecting direct pairs consisting of an $i\to0$ first-passage path of length $r$ and a $j\to0$ first-passage path of length $s$.
\end{lemma}

\begin{proof}
Take a crossed pair $(\omega,\omega')$ contributing to $F_r(j)F_s(i)$.  At time $0$, $\omega$ lies above $\omega'$, while at time $r$, $\omega_r=0$ and $\omega'_r>0$.  Since both paths live on the same parity lattice and move by nearest-neighbour steps, they intersect at some common time before $r$.  At the first such intersection, exchange the two tails.  The result is an intersecting direct pair contributing to $F_r(i)F_s(j)$.  The map is an involution, preserves first-passage admissibility, and preserves the product weight because the combined multiset of left and right steps is unchanged.  Hence crossed pairs cancel exactly against intersecting direct pairs, leaving only the nonintersecting direct pairs in the difference.
\end{proof}

\begin{remark}[Inhomogeneous chains and an algebraic proof]\label[remark]{rem:inhom}
The proof only exchanges tails at a common state, so the multiset of (state, step) transitions of the pair is preserved.  Hence \cref{lem:pathswitch} holds verbatim for nearest-neighbour chains with site-dependent probabilities $p_x,q_x$, with path weights given by products of transition probabilities.  This is used in \cref{sec:bd}.  An algebraic alternative: the two-step killed matrix $\mathsf Q$ restricted to one parity class is tridiagonal with $\mathsf Q(u,u)\mathsf Q(u{+}2,u{+}2)\ge \mathsf Q(u,u{+}2)\mathsf Q(u{+}2,u)$, hence $\TP$ (for a nonnegative tridiagonal matrix the adjacent $2\times2$ minors are the only ones that can be negative); lower exits occur only from the lowest state $i_0$ of the class, so $i\mapsto F_\cdot(i)$ is a column of $\mathsf Q^{k}$ times a constant, and Cauchy--Binet applied to $\mathsf Q^{k}[e_{i_0},\mathsf Q^{d}e_{i_0}]$, which uses only $2\times2$ minors, gives \eqref{eq:TP2}.
\end{remark}

For the symmetric path counts, reflection gives
\[
N_+(n;z)=N_-(n;2m-z).
\]
Taking $i=z$, $j=2m-z$, $r=n$, $s=n+2$ in \cref{lem:pathswitch} yields
\begin{equation}
N_-(n;z)N_+(n+2;z)-N_+(n;z)N_-(n+2;z)\ge0,
\label{eq:countdet}
\end{equation}
so $N_+/N_-$ is nondecreasing.  By \eqref{eq:biascancel}, $\eta_n/\nu_n$ is nondecreasing and $\varphi_n=1/(1+\eta_n/\nu_n)$ is nonincreasing.

\subsection{Strictness and the exceptional geometry}
Let $z<m$, put $j=2m-z$ and $d=m-z\ge1$.  Every post-opening time can be written
\[
n=j+2r,\qquad r\ge0.
\]
For $a\ge6$, the following two first-passage paths form a nonintersecting direct pair:
\begin{align*}
\omega:&\quad z,z-1,\ldots,1,
\underbrace{2,1,2,1,\ldots,2,1}_{d+r\text{ loops}},0,
\\
\omega':&\quad j,j-1,\ldots,3,
\underbrace{4,3,4,3,\ldots,4,3}_{1+r\text{ loops}},2,1,0.
\end{align*}
Their lengths are $n$ and $n+2$.  Until $\omega$ reaches the set $\{1,2\}$ the vertical separation is at least two; afterwards $\omega'$ stays at least at level $3$ until its final descent, while $\omega$ stays at levels $1,2$ before hitting $0$.  Thus the determinant in \eqref{eq:countdet} is strictly positive.

When $a=4,z=1$, any surviving path from $1$ begins $1\to2$ and any surviving path from $3$ begins $3\to2$.  Every direct pair therefore intersects at time one, so the determinant vanishes in every two-channel comparison.  This proves the strictness classification in \cref{thm:discreteMain}.  The limit \eqref{eq:philimit} follows from the principal eigenmode of the killed path matrix, whose reflected boundary amplitudes have equal leading path-count magnitude.

\begin{corollary}[Likelihood-ratio form]\label[corollary]{cor:evenlr}
Let $a=2m$ and $z<m$.  Then $G_0\le_\lr G_a$, strictly after both channels open when $a\ge6$.  Consequently all five responses of \cref{thm:dictionary} hold, in the direction favouring the lower boundary; for $z>m$ they hold in the reflected direction, and at $z=m$, $G_0=G_a$.
\end{corollary}

\begin{figure}[t]
\centering
\begin{tikzpicture}
\begin{axis}[
width=0.82\textwidth,height=6.1cm,
xlabel={admissible absorption time $n$},ylabel={$\varphi_n=\Pp(J=0\mid T=n)$},
ymin=-0.03,ymax=1.03,xmin=1,xmax=13,
legend style={at={(0.98,0.98)},anchor=north east,draw=none,fill=white},
grid=major]
\addplot+[mark=*,thick] coordinates {(2,1) (4,0.372) (6,0.270) (8,0.242) (10,0.233) (12,0.230)};
\addlegendentry{$a=6,p=0.6,z=2$}
\addplot+[mark=square*,thick] coordinates {(3,0.2285714) (5,0.2285714) (7,0.2285714) (9,0.2285714) (11,0.2285714) (13,0.2285714)};
\addlegendentry{$z=3$}
\addplot+[mark=triangle*,thick] coordinates {(2,0) (4,0.129) (6,0.192) (8,0.216) (10,0.224) (12,0.227)};
\addlegendentry{$z=4$}
\end{axis}
\end{tikzpicture}
\caption{Representative exact-time conditional lower-exit probabilities.  The lower-side initial point decreases toward the midpoint value, the midpoint is exactly neutral, and the reflected point increases.}
\label{fig:phiExample}
\end{figure}
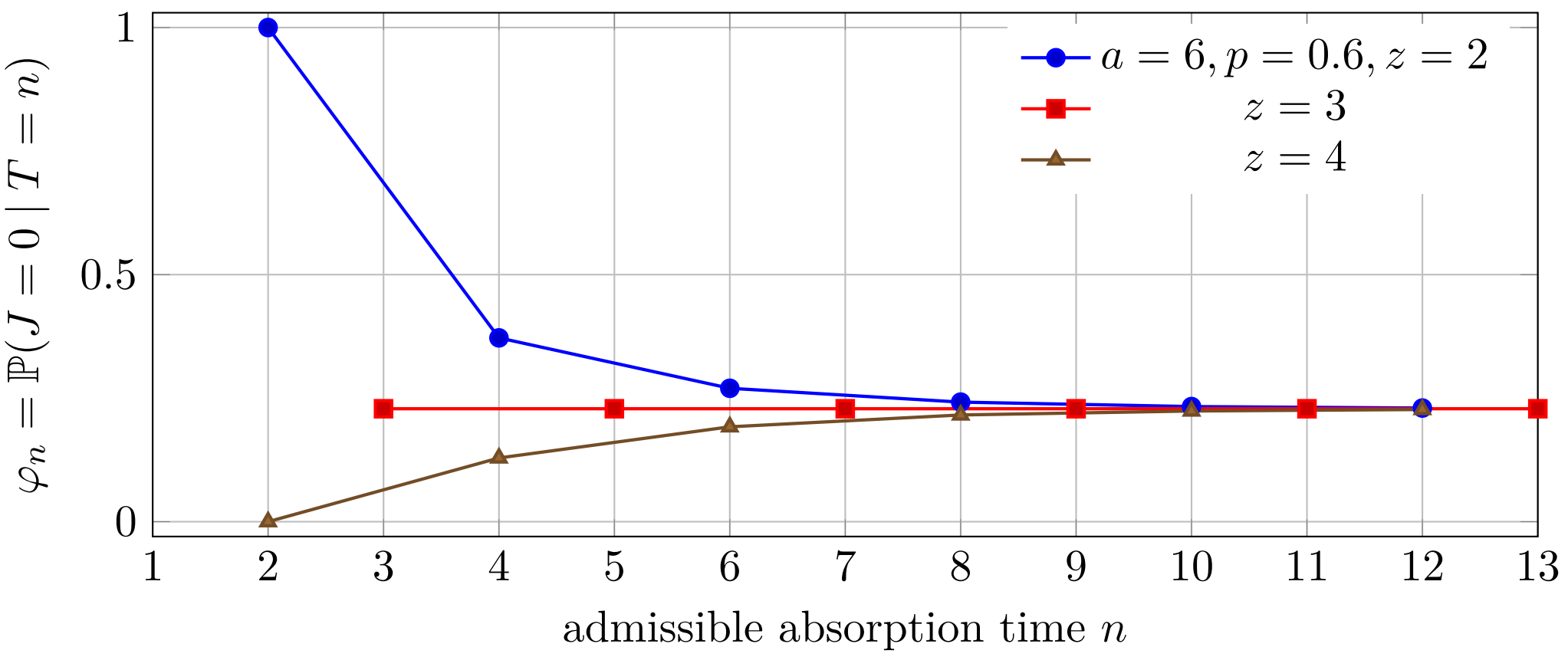

\section{Odd intervals: shifted likelihood ratios}\label{sec:odd}

If $a$ is odd, the two exit channels occur on opposite parities:
\[
T\equiv z\pmod2\quad\text{on }\{J=0\},
\qquad
T\equiv a-z\pmod2\quad\text{on }\{J=a\}.
\]
Hence $J$ is $\sigma(T)$-measurable and
\begin{equation}
H(J\mid T)=0,\qquad I(J;T)=H(J)
\label{eq:oddinfo}
\end{equation}
whenever both exits have nonzero probability.  The supports of $G_0$ and $G_a$ are interlaced, so $G_0\le_\lr G_a$ is impossible, and there is no universally reset-neutral site.  However, this does \emph{not} imply that every nontrivial reset mechanism changes the ruin probability.  For example, for the symmetric walk on $\{0,1,2,3\}$ started at $1$, the reset-free lower-exit probability is $2/3$.  A deterministic two-step cycle has one-cycle lower and upper absorption probabilities $1/2$ and $1/4$, giving the same renewal value $(1/2)/(3/4)=2/3$.  Universal neutrality is impossible because time determines side, although accidental equality may occur for individual reset laws.

The point of this section is that the orders which tolerate a one-step shift nevertheless survive.

\begin{lemma}[Tilt, reflection, first step]\label[lemma]{lem:tiltgen}
For every $a$, with $\rho=\sqrt{pq}$,
\[
\nu_n=\Big(\frac qp\Big)^{z/2}N_-(n;z)\,\rho^n,\qquad
\eta_n=\Big(\frac pq\Big)^{(a-z)/2}N_+(n;z)\,\rho^n .
\]
Moreover $N_+(n;z)=N_-(n;a-z)$, and for $n\ge2$,
$N_-(n;z)=N_-(n-1;z+1)+\ind_{\{z\ge2\}}N_-(n-1;z-1)$, with the convention $N_-(\cdot\,;a)=0$.
\end{lemma}

Thus the bias enters the exit pair only through the common tilt $\rho^n$ and constants; likelihood-ratio statements, shifted or not, are tilt invariant and reduce to path counts.  For even $a$ this is \eqref{eq:biascancel}.

\begin{lemma}[Shifted likelihood ratios]\label[lemma]{lem:shiftlr}
Let $a$ be odd and $z<a/2$.  For $k\equiv z\pmod 2$, $k\ge z$, put
\[
u_k=\frac{\eta_{k+1}}{\nu_k},\qquad v_k=\frac{\eta_{k-1}}{\nu_k}.
\]
Both sequences are nondecreasing along $k\equiv z\pmod2$.
\end{lemma}

\begin{proof}
By \cref{lem:tiltgen} it suffices to treat path counts; note $\nu_k>0$ for all such $k$.

\emph{$u_k$.} By reflection and the first step from $a-z$,
$N_+(k+1;z)=N_-(k+1;a-z)=\sum_{y\in Y}N_-(k;y)$ with $Y\subseteq\{a-z-1,\,a-z+1\}$.
Since $a$ is odd, every $y\in Y$ has the parity of $z$ and $y\ge a-z-1\ge z$.  By \cref{lem:pathswitch} (equality if $y=z$), $N_-(k;y)/N_-(k;z)$ is nondecreasing, and so is the sum.

\emph{$v_k$.} If $k=z=1$ then $v_1=0$; otherwise put $t=k-1\ge1$.  By the first step from $z$, $N_-(t+1;z)=\sum_{x\in X}N_-(t;x)$ with $X\subseteq\{z-1,\,z+1\}$, and $N_+(t;z)=N_-(t;a-z)$.  All of $x\in X$ and $y=a-z$ have the parity of $z+1$, and $x\le z+1\le y$.  \Cref{lem:pathswitch} gives $N_-(t;x)N_-(t+2;y)\ge N_-(t;y)N_-(t+2;x)$; summing over $x$,
\[
N_+(t+2;z)\,N_-(t+1;z)\ \ge\ N_+(t;z)\,N_-(t+3;z),
\]
which is $v_{k+2}\ge v_k$.
\end{proof}

\begin{theorem}[Odd intervals]\label{thm:odd}
Let $a$ be odd and $z<a/2$.
\begin{enumerate}[label=\textup{(\alph*)}]
\item $\Law(T-1\mid J=a)\ge_\lr G_0$.  Hence $G_0\le_\st G_a$ and $G_0\le_\pgf G_a$: every renewal reset law weakly increases the ruin probability, and under geometric resetting $Q_z(s)$ is nonincreasing in $s$.
\item Under deadlines, $Q_z(N+2)\le Q_z(N)$ for every $N\ge z$: the ruin probability is monotone along each parity class of $N$.
\item $G_0\le_\rh G_a$ fails in general, so $N\mapsto Q_z(N)$ need not be monotone.
\end{enumerate}
The inequalities reverse for $z>a/2$.
\end{theorem}

\begin{proof}
(a) is the monotonicity of $u_k$ together with \cref{rem:shift}.

(b) Let $B_N=\sum_{n\le N}\nu_n$ and $A_N=\sum_{n\le N}\eta_n$, so $Q_z(N)/(1-Q_z(N))=B_N/A_N$.  Passing from $N$ to $N+2$ adds one term to each sum; below $k$ runs over $k\equiv z$.
If $N\equiv z$, the new terms are $\eta_{N+1}=v_{N+2}\nu_{N+2}$ and $\nu_{N+2}$, while $A_N=\sum_{k\le N}v_k\nu_k\le v_{N+2}B_N$.
If $N\equiv z+1$, the new terms are $\eta_{N+2}=u_{N+1}\nu_{N+1}$ and $\nu_{N+1}$, while $A_N=\sum_{k\le N-1}u_k\nu_k\le u_{N+1}B_N$.
In both cases $A_N/B_N$ is at most the ratio of the added terms, so $A_{N+2}/B_{N+2}\ge A_N/B_N$ (mediant inequality), i.e.\ $Q_z(N+2)\le Q_z(N)$.

(c) $a=7$, $z=3$, $p=q$: $Q_z(4)=2/3<Q_z(5)=7/9$.
\end{proof}

\begin{theorem}[Hazard-rate ordered one-parity clocks]\label{thm:oneparity}
Let $a$ be odd and $z<a/2$.  Let $R_1,R_2$ be reset clocks supported on the same parity class of $\{z,z+1,\dots\}\cup\{\infty\}$, with survival sequences $w_1,w_2$.  If $w_1/w_2$ is nonincreasing, then
\[
Q_z(R_1)\ge Q_z(R_2).
\]
In particular the conclusion holds when $R_1\le_\lr R_2$.  The inequalities reverse for $z>a/2$.
\end{theorem}

\begin{proof}
Write the common support as $N_0<N_1<\cdots$ with $N_{k+1}-N_k=2$, and include $N_\infty=\infty$.  Compress time to this parity lattice.  Define conditional distribution functions
\[
\widetilde F_0(k)=\Pp(T\le N_k\mid J=0),\qquad
\widetilde F_a(k)=\Pp(T\le N_k\mid J=a).
\]
By \cref{thm:odd}(b), $Q_z(N_{k+1})\le Q_z(N_k)$, and the deadline-odds identity \eqref{eq:deadlineodds} is therefore equivalent to
$\widetilde F_a(k)/\widetilde F_0(k)$ being nondecreasing.  Thus the two compressed exit-time laws are ordered in reversed hazard rate.  A clock supported on the chosen parity class is exactly a clock on the compressed lattice, and the assumption that $w_1/w_2$ is nonincreasing is unchanged by the compression.  Applying the exact duality theorem \cref{thm:exactdual} on the compressed lattice gives the result.  Finally $\lr$ order of the clock laws implies the required hazard-rate order.
\end{proof}

\section{Reset ordering, attainable ranges, and the unified statement}\label{sec:resetorder}

\subsection{Survival-ratio ordering}
Let two reset clocks have survival sequences $w_1,w_2$.  Suppose
\begin{equation}
\frac{w_1(n)}{w_2(n)}\quad\text{is nonincreasing in }n.
\label{eq:survivalratio}
\end{equation}
Relative to clock 2, clock 1 gives progressively more weight to early absorption.  Equation \eqref{eq:renewalaverage} and monotone likelihood-ratio comparison give the following.

\begin{theorem}[Renewal reset-order theorem]\label{thm:resetorder}
For even $a=2m$, under \eqref{eq:survivalratio},
\[
Q_z(w_1)\ge Q_z(w_2),\qquad z<m,
\]
with the reverse inequality for $z>m$.  At $z=m$ equality holds for every pair of reset laws and every reset intensity.
\end{theorem}

\begin{proof}
For $z<m$, \cref{cor:evenlr} gives $G_0\le_\lr G_a$, and \cref{thm:dictionary}(v) gives the inequality; the case $z>m$ is its reflection.  At $z=m$, $G_0=G_a$ by \cref{cor:evenlr}, and \cref{thm:dictionary}(i) gives equality.
\end{proof}

For deterministic reset deadlines, increasing the number of permitted steps means weaker resetting.  The theorem therefore says that stronger resetting favours the nearer boundary.  For geometric survival $w_n=s^n$, decreasing $s$ increases reset intensity and moves the lower-exit probability upward for $z<m$.

\subsection{Correct attainable range}
A reset survival sequence is nonincreasing.  Thus a genuine reset mechanism cannot overweight arbitrarily late absorption relative to the no-reset law.  This removes the late-time limit $\varphi_*$ as an attainable extreme away from the midpoint.

\begin{corollary}[Sharp reset range]\label[corollary]{cor:range}
Include the degenerate no-reset law $w_n\equiv1$.  For even $a=2m$,
\begin{align*}
 z<m:&\qquad Q_z(w)\in[Q_z^{(0)},1],\\
 z=m:&\qquad Q_z(w)=\frac{q^m}{p^m+q^m},\\
 z>m:&\qquad Q_z(w)\in[0,Q_z^{(0)}],
\end{align*}
where $Q_z^{(0)}$ is the reset-free lower-exit probability.  These intervals are sharp.
\end{corollary}

\begin{proof}
The ordering against $w\equiv1$ gives the inner endpoint.  If $z<m$, choose a deterministic deadline after the lower boundary is reachable but before the upper boundary is reachable; then every successful cycle exits at $0$, so the eventual lower-exit probability is $1$.  The reflected construction gives $0$ for $z>m$.  The geometric family varies continuously between the no-reset value and the near-boundary extreme, so every intermediate value is attained.
\end{proof}

\begin{remark}
By \cref{thm:odd}(a) and the same deadline construction, \cref{cor:range} holds verbatim for odd $a$ with $z<a/2$ (and in reflected form for $z>a/2$).
\end{remark}

\subsection{The unified statement}

\begin{theorem}[Parity-corrected likelihood-ratio order]\label{thm:unified}
For the nearest-neighbour walk on $\{0,\dots,a\}$ with $p,q>0$ and $z<a/2$,
\[
\Law\big(T-(a\bmod 2)\mid J=a\big)\ \ge_\lr\ \Law(T\mid J=0).
\]
\end{theorem}

\begin{center}
\small
\begin{tabular}{@{}p{0.43\textwidth}p{0.19\textwidth}p{0.30\textwidth}@{}}
\toprule
Reset response ($z<a/2$) & $a$ even & $a$ odd\\
\midrule
Universal sign, $Q_z\ge\pi_0$ ($\st$) & yes & yes\\
Geometric monotone ($\pgf$) & yes & yes\\
Deadline monotone ($\rh$) & yes & no; yes along each parity of $N$\\
Hazard-rate ordered clocks ($\rh$, \cref{thm:exactdual}) & yes & no in general; yes on one parity\\
Arbitrary gated weights ($\lr$) & yes & no; shifted by one step: yes\\
Universal neutrality & only at $z=m$ & never (accidental equality possible)\\
\bottomrule
\end{tabular}
\end{center}

\section{General birth--death chains: the switching site}\label{sec:bd}

Let $X$ be a nearest-neighbour chain on $\{0,\dots,a\}$, absorbed at $0$ and $a$, with $\Pp(x\to x+1)=p_x$, $\Pp(x\to x-1)=q_x$, $p_x,q_x>0$, $p_x+q_x=1$ for $0<x<a$.  Let $P$ be the killed one-step matrix on the interior.  The chain is reversible with respect to $\pi(x)\propto\prod_{0<y<x}p_y/q_{y+1}$.  By \cref{rem:inhom}, \cref{lem:pathswitch} holds for $X$.  Throughout this section $a$ is even, so both exits from $z$ occur on the parity $n\equiv z$; the odd case is analogous with the shifted ratios of \cref{lem:shiftlr}.  Write $\theta_z(n)=\eta_n/\nu_n$ for the exact-time odds of the far side.

\begin{lemma}[Arrival representation]\label[lemma]{lem:arrival}
There is $c_\ast>0$, independent of $z$ and $n$, such that
\[
\theta_z(n)=c_\ast\,\frac{P^{\,n-1}(a-1,z)}{P^{\,n-1}(1,z)} .
\]
\end{lemma}

\begin{proof}
One has
\[
\nu_n=q_1P^{n-1}(z,1),\qquad
\eta_n=p_{a-1}P^{n-1}(z,a-1).
\]
Reversibility gives
\[
P^k(z,i)=\frac{\pi(i)}{\pi(z)}P^k(i,z),\qquad
c_\ast=\frac{p_{a-1}\pi(a-1)}{q_1\pi(1)},
\]
which yields the stated representation.
\end{proof}

The exit odds from $z$ are therefore a ratio of probabilities of \emph{arriving} at $z$ from the two extreme interior states.

\begin{theorem}[Monotone switch]\label{thm:switch}
For admissible $n$ with both channels open, let $D(n,z)=\log\theta_z(n+2)-\log\theta_z(n)$.  For every $n$, $z\mapsto D(n,z)$ is nonincreasing on each parity class of starting sites.  Consequently, at each time the sites whose odds increase form an initial segment of the class and those whose odds decrease form a final segment; in particular, the sites with globally nondecreasing odds form an initial segment and those with globally nonincreasing odds a final segment.
\end{theorem}

\begin{proof}
Let $x_1<x_2$ have the same parity.  \Cref{lem:pathswitch} for lower exits, $\nu_n(x)=q_1P^{n-1}(x,1)$, together with reversibility (the factors $\pi(x)$ cancel in ratios), shows that $P^k(1,x_2)/P^k(1,x_1)$ is nondecreasing along admissible $k$.  The same lemma for upper exits, i.e.\ for the reversed order of states, shows that $P^k(a-1,x_2)/P^k(a-1,x_1)$ is nonincreasing.  By \cref{lem:arrival},
\[
\frac{\theta_{x_2}(n)}{\theta_{x_1}(n)}
=\frac{P^{n-1}(a-1,x_2)/P^{n-1}(a-1,x_1)}{P^{n-1}(1,x_2)/P^{n-1}(1,x_1)}
\]
is nonincreasing in $n$, i.e.\ $D(n,x_2)\le D(n,x_1)$.
\end{proof}

\begin{corollary}[Antisymmetry; a reflection-free proof of the even theorem]\label[corollary]{cor:antisym}
Suppose an order-reversing involution $\sigma$ of a parity class satisfies $D(n,\sigma z)=-D(n,z)$.  Then $D(n,z)\ge0$ whenever $z<\sigma z$.  For the nearest-neighbour walk with $a$ even, $\sigma z=a-z$ preserves each parity class.  By \cref{lem:tiltgen}, $D$ does not depend on $p$, and for $p=q$ symmetry gives $\theta_{a-z}=1/\theta_z$.  This recovers the monotonicity part of \cref{thm:discreteMain} without the reflection step.
\end{corollary}

\begin{proof}
If $z<\sigma z$, \cref{thm:switch} gives $D(n,z)\ge D(n,\sigma z)=-D(n,z)$.
\end{proof}

\begin{proposition}[Asymptotic switching site]\label[proposition]{prop:node}
Let $C$ be the parity class of $z$, of size $r\ge2$, and let $\mathsf Q$ be the two-step killed matrix restricted to $C$. Assume its two largest eigenvalues satisfy $\mu_1>\mu_2>0$, and let $\psi_1,\psi_2$ be the corresponding eigenvectors of its symmetrisation, with $\psi_1>0$ and $\psi_2(\min C)>0$.  If $\psi_2(z)\neq 0$, then for all sufficiently large admissible $n$
\[
\sgn D(n,z)=\sgn\psi_2(z).
\]
Consequently the eventual sign follows the single sign change of $\psi_2$ on sites where $\psi_2$ is nonzero.  If $\mu_2=0$ or $\psi_2(z)=0$, the displayed leading-term argument does not determine an eventual nonzero sign; the $a=4$ exceptional geometry has $D\equiv0$ after opening.
\end{proposition}

\begin{proof}
$\mathsf Q$ is tridiagonal with positive off-diagonal entries, so its spectrum is simple; it is a product of the two bidiagonal one-step blocks between the parity classes, hence totally nonnegative and $\mu_j\ge0$.  By \cref{lem:arrival}, $\theta_z\propto \mathsf Q^k(i_t,z)/\mathsf Q^k(i_0,z)$ with $i_0=\min C$, $i_t=\max C$ and $k$ the block time.  If $z$ is odd this is immediate, since $1,a-1\in C$.  If $z$ is even, the first step from $1$ (resp.\ $a-1$) is forced to $2$ (resp.\ $a-2$), so $P^{n-1}(1,z)=p_1P^{n-2}(2,z)$ and $P^{n-1}(a-1,z)=q_{a-1}P^{n-2}(a-2,z)$, and the constant factors cancel in $D$.  By the oscillation theorem for Jacobi matrices \cite{GantmacherKrein}, $\psi_2$ changes sign exactly once and does not vanish at the end points, so $\psi_2(i_0)>0>\psi_2(i_t)$.  With $\varepsilon=\mu_2/\mu_1<1$ and $\beta_i=\psi_2(i)/\psi_1(i)$, the spectral expansion gives
\[
D=(\varepsilon-1)\,\varepsilon^k\,\beta_z\,(\beta_{i_t}-\beta_{i_0})+o(\varepsilon^k),
\]
and both $\varepsilon-1$ and $\beta_{i_t}-\beta_{i_0}$ are negative.
\end{proof}

\begin{remark}[Degenerate parity blocks]
If a parity block has a single state, its exact-time exit odds are constant wherever both channels are open.  A second mode is then undefined.  If the second eigenvalue vanishes, the leading asymptotic correction in \cref{prop:node} vanishes: for example, in the symmetric walk on $\{0,1,2,3,4\}$, the odd-class two-step matrix is $\frac14\left(\begin{smallmatrix}1&1\\1&1\end{smallmatrix}\right)$, with eigenvalues $1/2$ and $0$.  Starting at $z=1$, the post-opening odds are constant and $D(n,1)=0$ for all odd $n\ge3$ despite $\psi_2(1)\ne0$.
\end{remark}

\begin{proposition}[Failure of the single-turn property]\label[proposition]{prop:singlefail}
The map $n\mapsto D(n,z)$ can change sign more than once for an inhomogeneous birth--death chain.  In particular, take $a=12$, $z=5$, and
\[
(p_1,\ldots,p_{11})=
\left(\frac35,\frac15,\frac{49}{50},\frac1{50},\frac12,\frac12,
\frac{19}{20},\frac12,\frac25,\frac1{10},\frac12\right).
\]
Then $D(7)>0$, $D(9)>0$, $D(11)<0$, $D(13)<0$, and $D(15)>0$.  Thus the exact-time odds increase, then decrease, then increase again.
\end{proposition}

\begin{proof}
Since the logarithm is increasing, the sign of $D(n,5)$ is the sign of
\[
\Delta_n:=\eta_{n+2}\nu_n-\eta_n\nu_{n+2}.
\]
The one-step recurrence with the rational transition probabilities above gives, exactly,
\[
\begin{aligned}
\operatorname{numerator}(\Delta_7)&=4901373323,\\
\operatorname{numerator}(\Delta_9)&=467939702519611,\\
\operatorname{numerator}(\Delta_{11})&=-24101569376877744983,\\
\operatorname{numerator}(\Delta_{13})&=-13053045806915174743881261,\\
\operatorname{numerator}(\Delta_{15})&=4509048914570924507839832768453.
\end{aligned}
\]
All omitted denominators are positive.  The displayed sign pattern is therefore exact.
\end{proof}

\begin{proposition}[A spectral bound on temporal switches]\label[proposition]{prop:switchbound}
Let $C$ be the parity class of $z$ with $r=|C|\ge2$, $i_0=\min C$, $i_t=\max C$, and let $\mu_1>\dots>\mu_r\ge0$ be the eigenvalues of the two-step killed matrix $\mathsf Q$ on $C$, with orthonormal eigenvectors $\psi_j$ of its symmetrisation.  For the two-channel block clock $k\ge1$ (two ordinary steps per block), $\sgn D=\sgn W(k)$, where
\[
W(k)=\sum_{i<j}c_{ij}\,(\mu_i\mu_j)^k,\qquad
c_{ij}=(\mu_i-\mu_j)\,\psi_i(z)\psi_j(z)\big(\psi_i(i_t)\psi_j(i_0)-\psi_j(i_t)\psi_i(i_0)\big).
\]
Consequently, on the two-channel block clock, the number of strict sign reversals of $D(n,z)$ (zeros omitted) is at most the sign variation of the nonzero coefficients with positive bases $\mu_i\mu_j$, after merging equal bases and ordering them decreasingly. If $W\not\equiv0$, this is at most $\binom r2-1$; if $W\equiv0$, there are no strict reversals. The case $r=1$ also has no reversals. Under the additional hypotheses of \cref{prop:node}, the eventual sign is $\sgn\psi_2(z)$.
\end{proposition}

\begin{proof}
By \cref{lem:arrival} and the spectral decomposition of the symmetrisation, the lower and upper exits at block $k+1$ are positive multiples, independent of $k$, of $g(k)=\sum_j\psi_j(z)\psi_j(i_0)\mu_j^k$ and $h(k)=\sum_j\psi_j(z)\psi_j(i_t)\mu_j^k$.  Hence $\sgn D=\sgn\big(h(k+1)g(k)-h(k)g(k+1)\big)$, and expanding and pairing the terms $(i,j)$ and $(j,i)$ gives $W(k)$.  Terms with $\mu_i\mu_j=0$ vanish for $k\ge1$.  After merging equal bases, $t\mapsto W(t)$ is a real exponential sum $\sum_lC_l\lambda_l^t$ with distinct $\lambda_l>0$, so by the Descartes rule of signs for exponential sums \cite[Part~V]{PolyaSzego} its number of real zeros is at most the number of sign changes of $(C_l)$, ordered by $\lambda_l$.  A sign change of the sequence $W(k)$ forces a real zero between consecutive integers.
\end{proof}

\begin{remark}
For the chain of \cref{prop:singlefail} ($r=6$) the coefficient bound gives at most $5$ switches, against $\binom62-1=14$; exactly two reversals occur globally by \cref{cor:exacttwoswitch}.  Thus a site may switch direction several times, but only finitely often, by an amount controlled by the spectrum; whenever the hypotheses of \cref{prop:node} hold, the eventual sign is governed by $\psi_2$.  The gap between the certified exact count $2$ and the coefficient bound $5$ is part of the open problem in \cref{sec:open}.
\end{remark}

\begin{corollary}[Exactly two temporal reversals in the rational example]\label[corollary]{cor:exacttwoswitch}
For the chain in \cref{prop:singlefail}, with $z=5$, the exact-time log-odds increment satisfies
\[
D(7,5)>0,\quad D(9,5)>0,\quad
D(11,5)<0,\quad D(13,5)<0,
\qquad D(n,5)>0\quad\text{for every odd }n\ge15.
\]
Consequently it has \emph{exactly two} strict temporal sign reversals over its entire admissible two-channel time domain, at the transitions to $n=11$ and $n=15$.
\end{corollary}

\begin{proof}
The initial six signs at $n=7,9,11,13,15,17$ follow by direct rational powers of the six-state two-step killed matrix.  The remaining infinite tail is certified by rational interval arithmetic: the five positive eigenvalues are enclosed in disjoint intervals with exact polynomial sign changes at their endpoints, and the dominant coefficient in the exponential sum of \cref{prop:switchbound} strictly exceeds the combined absolute contribution of all other coefficients for every block index $k\ge9$.  The full explicit characteristic polynomial, cofactor coefficients, rational root brackets and tail inequality are given in \cref{app:exacttwoswitch}; the standalone exact replay script performs the same checks independently of floating-point eigensolvers.
\end{proof}

\begin{remark}[Minimality of the length]\label[remark]{rem:minlength}
The working manuscript \cite{NewtonReversals} of the first author reports that $a=12$ is the smallest even interval length admitting an inhomogeneous birth--death chain with two temporal reversals.
\end{remark}

\section{Response identities, susceptibility, and information structure}

\subsection{General covariance response}
Let $w_\theta(n)>0$ be a differentiable family of reset survival weights and define the tilted absorption-time law
\[
\varpi_\theta(n)=\frac{m_nw_\theta(n)}{\sum_km_kw_\theta(k)}.
\]
Then $Q_z(\theta)=\Ee_{\varpi_\theta}[\varphi_T]$.  Differentiating a normalized weighted mean yields
\begin{theorem}[Reset-response covariance identity]\label{thm:covresp}
\begin{equation}
\boxed{
\frac{\dd Q_z}{\dd\theta}
=\Cov_{\varpi_\theta}\!\left(
\varphi_T,\,\partial_\theta\log w_\theta(T)
\right).}
\label{eq:covresp}
\end{equation}
\end{theorem}

For the exponential/geometric parametrization $w_\theta(n)=e^{-\theta n}$,
\begin{equation}
\frac{\dd Q_z}{\dd\theta}
=-\Cov_{\varpi_\theta}(\varphi_T,T)
=-\Cov_{\varpi_\theta}(\ind_{\{J=0\}},T).
\notag
\end{equation}
Thus the infinitesimal reset susceptibility is exactly minus the exit-side/absorption-time covariance under the tilted cycle law.

\subsection{Exact symmetric-walk susceptibility}
For the unbiased walk $p=q=1/2$, at zero reset intensity,
\begin{equation}
\boxed{
-\Cov_z(\ind_{\{J=0\}},T)
=\frac{z(a-z)(a-2z)}{3a}.}
\label{eq:discretesusc}
\end{equation}
A direct proof follows from the standard second-difference equations for $\Ee_zT$ and $\Ee_z[T\ind_{\{J=0\}}]$; details are given in Appendix~B.  The sign changes once, at the midpoint, exactly as predicted by \cref{thm:discreteMain}.

\subsection{Information-theoretic characterization}
At an even midpoint,
\[
J\perp T,
\]
so
\begin{equation}
I(J;T)=0.
\notag
\end{equation}
For odd intervals, \eqref{eq:oddinfo} gives the opposite extreme.  The parity transition can therefore be summarized as
\[
\boxed{
\text{even midpoint: }I(J;T)=0,
\qquad
\text{odd interval: }I(J;T)=H(J).}
\]
Reset neutrality is consequently not only a symmetry of a renewal formula: it is a precise decoupling of temporal and spatial exit information.

\section{Drift-diffusion: strict continuum monotonicity and exact reset response}\label{sec:BM}

\subsection{Boundary first-passage densities}
Consider
\begin{equation}
\dd X_t=\mu\,\dd t+\sigma\,\dd B_t,
\qquad 0<X_t<a,
\notag
\end{equation}
with absorbing endpoints and $X_0=x$.  Let $g_0(t;x)$ and $g_a(t;x)$ be the boundary-specific first-passage densities.  The killed spectral kernel gives
\begin{align}
 g_0(t;x)
 &=\frac{\sigma^2\pi}{a^2}
 e^{-\mu x/\sigma^2-\mu^2t/(2\sigma^2)}
 \sum_{n\ge1}n\sin\frac{n\pi x}{a}
 e^{-\sigma^2n^2\pi^2t/(2a^2)},
\notag\\
 g_a(t;x)
 &=\frac{\sigma^2\pi}{a^2}
 e^{\mu(a-x)/\sigma^2-\mu^2t/(2\sigma^2)}
 \sum_{n\ge1}(-1)^{n+1}n\sin\frac{n\pi x}{a}
 e^{-\sigma^2n^2\pi^2t/(2a^2)}.
\notag
\end{align}
Hence the drift factors out of the exact-time odds:
\begin{equation}
\frac{g_0(t;x)}{g_a(t;x)}
=e^{-\mu a/\sigma^2}
\frac{h_0(t;x)}{h_0(t;a-x)},
\label{eq:contbiascancel}
\end{equation}
where $h_0$ is the zero-drift lower first-passage density.

Karlin and McGregor proved total positivity for classical one-dimensional diffusion kernels \cite{KarlinMcGregor1960}; Karlin's first-passage theorem \cite[Theorem~5.1(i)]{Karlin1964} allows defective first-passage distributions.  Here the relevant density is specifically
$h_0(t;x)\,\dd t=\Pp_x(T_0\in\dd t,\,T_0<T_a)$, with killing at \emph{both} endpoints.  We avoid any ambiguity about the applicability of a one-boundary theorem by deriving its strict order directly from the discrete first-passage inequality through an explicit spectral local limit.

\begin{lemma}[Parity-resolved killed-walk local limit]\label[lemma]{lem:killedlocallimit}
Fix $a,\sigma>0$, $x\in(0,a)$ and $t>0$.  For every integer $N\ge2$ let $\delta_N=a^2/(\sigma^2N^2)$, and consider the simple symmetric nearest-neighbour walk on $\{0,\dots,N\}$ absorbed at its endpoints.  Choose interior starting sites $i_N$ such that $ai_N/N\to x$ and admissible positive integers $n_N$ with $n_N\equiv i_N\pmod2$ and $n_N\delta_N\to t$.  Then the probability $F_{n_N}^{(N)}(i_N)$ of first exit through zero at step $n_N$ satisfies
\begin{equation}\label{eq:locallimit}
\frac{F_{n_N}^{(N)}(i_N)}{2\delta_N}\ \longrightarrow\
h_0(t;x)=\frac{\sigma^2\pi}{a^2}\sum_{j\ge1} j\sin\!\left(\frac{j\pi x}{a}\right)
\exp\!\left(-\frac{\sigma^2j^2\pi^2t}{2a^2}\right).
\end{equation}
The same convergence holds simultaneously for finitely many distinct starting points and times, while keeping all starting sites on one parity class.
\end{lemma}
\begin{proof}
Write $P_N$ for the interior killed-walk matrix.  Its orthogonal sine eigenbasis gives, exactly,
\begin{equation}\notag
F_n^{(N)}(i)=\tfrac12 P_N^{n-1}(i,1)
=\frac1N\sum_{j=1}^{N-1}
\sin\!\left(\frac{j\pi i}{N}\right)
\sin\!\left(\frac{j\pi}{N}\right)
\cos^{n-1}\!\left(\frac{j\pi}{N}\right).
\end{equation}
For each fixed $j$, the low-frequency term, divided by $2\delta_N$, tends to half the $j$th summand in \eqref{eq:locallimit}.  The complementary mode $N-j$ has the same sign and limiting contribution because $n\equiv i\pmod2$; the two contributions supply exactly the factor $2$ in the time cell $2\delta_N$.  Indeed, for $j\le N/2$, the bound
$|\sin(j\pi/N)|\le\pi j/N$, together with
$|\cos(j\pi/N)|^{n_N-1}\le \exp(-c_tj^2)$ for all sufficiently large $N$ and some $c_t>0$, gives a summable bound $Cj e^{-c_tj^2}$ on each rescaled term.  Symmetry $j\leftrightarrow N-j$ and dominated convergence justify the limit of the entire sum.  Choosing common-parity approximating sites and times is always possible by changing each lattice index by at most one.
\end{proof}

\begin{proposition}[Strict first-passage order for the two-sided killed kernel]\label[proposition]{prop:killedTP2}
For Brownian motion of volatility $\sigma>0$ absorbed on leaving $(0,a)$, the lower-exit density satisfies
\begin{equation}\label{eq:killedfirstpassTP2}
h_0(t;x)h_0(s;y)>h_0(s;x)h_0(t;y),
\qquad 0<x<y<a,\quad 0<t<s.
\end{equation}
\end{proposition}
\begin{proof}
For $0<x<y<a$ choose lattice starting sites $i_N<j_N$ of the same parity and approximating $x,y$.  Similarly choose $r_N<s_N$ of that parity with $r_N\delta_N\to t$ and $s_N\delta_N\to s$.  The exact lattice path-switching lemma \ref{lem:pathswitch} gives
$F_{r_N}^{(N)}(i_N)F_{s_N}^{(N)}(j_N)-F_{s_N}^{(N)}(i_N)F_{r_N}^{(N)}(j_N)\ge0$.
Divide by $(2\delta_N)^2$ and apply \cref{lem:killedlocallimit}.  This proves the weak version of \eqref{eq:killedfirstpassTP2} with the \emph{two-sided} killing already incorporated.  Because the densities are strictly positive and real analytic in $t>0$, the ratio $h_0(t;x)/h_0(t;y)$ is positive, analytic and nonincreasing.  Its short-time boundary-flux asymptotics give
\[
h_0(t;u)\sim \frac{u}{\sqrt{2\pi\sigma^2t^3}}
\exp\!\left(-\frac{u^2}{2\sigma^2t}\right)
\qquad (t\downarrow0),\quad 0<u<a,
\]
so the ratio diverges to $+\infty$ as $t\downarrow0$ when $x<y$.  The principal Dirichlet eigenfunction gives the finite positive limit
$\sin(\pi x/a)/\sin(\pi y/a)$ as $t\to\infty$.  An analytic nonincreasing function taking the same value at two distinct positive times would be constant on the intervening interval, hence everywhere, contradicting these limits.  Thus it is strictly decreasing, proving \eqref{eq:killedfirstpassTP2}.
\end{proof}

For $x<a/2$, set $y=a-x$ in \eqref{eq:killedfirstpassTP2}.  Reflection identifies $h_0(t;y)$ with the zero-drift upper-exit density at $x$.  Consequently the drift-independent ratio in \eqref{eq:contbiascancel} is strictly decreasing.  This establishes the continuum response without an appeal to unverified killed-semigroup hypotheses of the cited classical theorem.
\begin{theorem}[Continuum exact-time monotonicity]\label{thm:continuumphi}
Define
\[
\Phi_x(t)=\Pp_x(J=0\mid\tau=t)
=\frac{g_0(t;x)}{g_0(t;x)+g_a(t;x)}.
\]
Then $\Phi_x(t)$ is strictly decreasing for $x<a/2$, constant for $x=a/2$, and strictly increasing for $x>a/2$.  Moreover
\begin{equation}
\lim_{t\to\infty}\Phi_x(t)
=\frac{1}{1+e^{\mu a/\sigma^2}}
\label{eq:continuumlimit}
\end{equation}
for every $x$.
\end{theorem}

\begin{proof}
By \eqref{eq:contbiascancel}, $\Phi_x(t)=\bigl(1+e^{\mu a/\sigma^2}h_0(t;a-x)/h_0(t;x)\bigr)^{-1}$.  For $x<a/2$, \cref{prop:killedTP2} with $y=a-x$ shows that $h_0(t;x)/h_0(t;a-x)$ is strictly decreasing, so $\Phi_x$ is strictly decreasing; $x>a/2$ is the reflection, and at $x=a/2$ the ratio is $1$.  As $t\to\infty$ the principal Dirichlet mode gives $h_0(t;a-x)/h_0(t;x)\to\sin(\pi(a-x)/a)/\sin(\pi x/a)=1$, which is \eqref{eq:continuumlimit}.
\end{proof}

Unlike the lattice walk, Brownian motion has no initial one-channel plateau: each boundary has positive first-passage density for every $t>0$.

\subsection{General renewal reset ordering}
For an independent continuous reset time with survival function $W(t)$,
\begin{equation}
Q_x(W)=
\frac{\int_0^\infty W(t)g_0(t;x)\dd t}
{\int_0^\infty W(t)[g_0(t;x)+g_a(t;x)]\dd t}.
\notag
\end{equation}
Thus the discrete weighted-average argument passes directly to the continuum.  If $W_1/W_2$ is nonincreasing, clock 1 is an early-time likelihood-ratio tilt and \cref{thm:continuumphi} gives
\[
Q_x(W_1)\ge Q_x(W_2),\quad x<a/2,
\]
with the reverse inequality for $x>a/2$.

At the midpoint,
\begin{equation}
Q_{a/2}(W)=\frac{1}{1+e^{\mu a/\sigma^2}}
\notag
\end{equation}
for every independent renewal reset law.

\subsection{Poisson resetting and exact susceptibility}
For Poisson reset rate $\lambda$, $W(t)=e^{-\lambda t}$.  Put
\[
\kappa_\lambda=\frac{\sqrt{\mu^2+2\sigma^2\lambda}}{\sigma^2}.
\]
The discounted boundary probabilities are
\begin{align*}
A_0(\lambda)&=e^{-\mu x/\sigma^2}
\frac{\sinh(\kappa_\lambda(a-x))}{\sinh(\kappa_\lambda a)},\\
A_a(\lambda)&=e^{\mu(a-x)/\sigma^2}
\frac{\sinh(\kappa_\lambda x)}{\sinh(\kappa_\lambda a)},
\end{align*}
and renewal gives $Q_x(\lambda)=A_0/(A_0+A_a)$.  This is the drifted analogue of the interval splitting formulas studied in \cite{PalPrasad2019}.

Differentiation yields
\begin{equation}
\boxed{
\frac{\dd Q_x}{\dd\lambda}
=\frac{Q_x(1-Q_x)}{\sigma^2\kappa_\lambda}
\left[(a-x)\coth(\kappa_\lambda(a-x))
-x\coth(\kappa_\lambda x)\right].}
\notag
\end{equation}
Since $s\mapsto s\coth(\kappa_\lambda s)$ is strictly increasing,
\[
\sgn Q_x'(\lambda)=\sgn(a-2x).
\]
There is therefore no interior optimum of the splitting probability as a function of reset rate: increasing reset intensity monotonically favours the nearer boundary.

For zero drift,
\begin{equation}
Q_x'(0)=\frac{x(a-x)(a-2x)}{3a\sigma^2}
=-\Cov_x(\ind_{\{J=0\}},\tau).
\label{eq:contsusc}
\end{equation}
Equation \eqref{eq:contsusc} is the diffusive scaling counterpart of \eqref{eq:discretesusc}.

\section{Concluding remarks and open problems}\label{sec:open}

\subsection{Exact certificates}
The proof of \cref{thm:discreteMain} is analytic.  Independently, an exact integer recurrence
\[
c_{n+1}(x)=c_n(x-1)+c_n(x+1)
\]
was used as a regression certificate.  For $a=4,6,\ldots,100$, $z<a/2$ and admissible $n<499$ there were $295{,}525$ determinant comparisons of the form \eqref{eq:countdet}; none were negative.  Of these, $275{,}925$ comparisons had the far-exit channel accessible by the \emph{later} time $n+2$.  Only $248$ determinants in that group vanished, all with $(a,z)=(4,1)$ and both channels already accessible at the earlier time.  These endpoint conventions matter to the reproducibility of the counts.  Further audits are listed in the appendix.

\subsection{Prior art}
Total positivity and first-passage results for Markov chains and diffusions are classical \cite{Karlin1964,KarlinMcGregor1960,Karlin1968}; stochastic orders and their implications are standard \cite{ShakedShanthikumar2007}; renewal resetting and interval splitting probabilities have a substantial literature \cite{EvansMajumdar2011,PalPrasad2019,Belan2018,EvansMajumdarSchehr2020}; universal orderings of a single completion time under restart were characterised in \cite{Vidmar2022}, and moment criteria for selected outcomes were given in \cite{PalEtAl2025}.  The contributions of this paper are the dictionary between reset responses and stochastic orders of the exit pair, its Pl\"ucker formulation and exact duality, the self-contained discrete monotonicity theorem with strictness classification, the parity-corrected theorem for odd intervals, the switching theorem for general birth--death chains, and the continuum counterpart for drifted Brownian motion.

\subsection{Open problems}
\begin{enumerate}
\item \emph{Temporal switching complexity.} \Cref{cor:exacttwoswitch} certifies a general inhomogeneous chain with exactly two temporal switches, at the smallest even length that admits them (\cref{rem:minlength}).  \Cref{prop:switchbound} bounds the strict sign reversals of $n\mapsto D(n,z)$ by those of a spectral coefficient sequence (at most $\binom r2-1$ when $r\ge2$).  Determine sharp bounds, and identify natural structural hypotheses on the environment $(p_x)$ under which a genuine single-turn theorem is restored.  Second-compound and variation-diminishing methods are natural candidates.
\item \emph{Duality beyond hazard rate.} \Cref{thm:exactdual} settles hazard-rate and likelihood-ratio ordered clocks, while \cref{prop:clockstfail} shows that ordinary stochastic order is too weak.  Open: the exact dual of other classical clock orders, including mean residual life and convex orders.
\item \emph{Optimal resetting at multiply switching sites.}  The failure of the single-turn property permits richer finite-intensity response curves than previously expected.  Classify when a birth--death chain can have one or several interior extrema as the reset intensity varies.  Drifted Brownian motion remains the contrasting no-interior-optimum case (\cref{sec:BM}).
\item \emph{Continuum switching.} The proof of \cref{thm:switch} uses only $\TP$ and reversibility; one-dimensional diffusions are reversible with respect to their speed measure, so a continuum switching theorem with switching point at the node of the second Dirichlet eigenfunction is expected.
\item \emph{Continuum dictionary for general diffusions.} \Cref{sec:BM} treats drifted Brownian motion.  For variable coefficients the likelihood-ratio order below the midpoint cannot be expected in general, since \cref{prop:singlefail} already exhibits multiple reversals in the discrete setting.  It should hold, with the full dictionary and its strictness classification, when the Liouville potential is symmetric about the intrinsic midpoint and the drift is arbitrary, by the Karlin--McGregor coincidence formula \cite{KarlinMcGregor1960}, the continuum counterpart of the path-switching lemma.  For general coefficients the natural statements are the continuum switching theorem of the previous item and its late-time node.
\end{enumerate}

The conclusion is that the direction of the reset response is a property of the exit pair alone: comparing two reset mechanisms is a pairing on $\Lambda^2$, and the classical stochastic orders are exactly the positivity conditions that this pairing requires.

\appendix
\section{Main theorem package}
\subsection*{A.1 Discrete theorem package}
Let $X_n$ be the absorbed nearest-neighbour walk on $\{0,\ldots,2m\}$ with $p,q>0$, started at $z$.  Let $T$ be absorption time, $J\in\{0,2m\}$ the exit side, and $\varphi_n=\Pp(J=0\mid T=n)$.  Let reset cycles be i.i.d., independent of the walk, and described by an admissible survival sequence $w_n=\Pp(R\ge n)$ with $\sum_nm_nw_n>0$.

The precise discrete claims are:
\begin{enumerate}
\item Renewal identity: $Q_z(w)=\sum\varphi_nm_nw_n/\sum m_nw_n$.
\item If $z<m$, then $\varphi_n=1$ before time $2m-z$, then $\varphi_n$ is strictly decreasing on the admissible parity for all $a\ge6$; $(a,z)=(4,1)$ is the unique strictness exception.  Reflection gives the $z>m$ statement.  At $z=m$, $\varphi_n=q^m/(p^m+q^m)$ for all admissible $n$.
\item If $w_1/w_2$ is nonincreasing, then $Q_z(w_1)\ge Q_z(w_2)$ for $z<m$, with the reverse inequality for $z>m$.
\item Universal renewal neutrality is equivalent to $J\perp T$.
\item For differentiable survival weights, \eqref{eq:covresp} holds.
\end{enumerate}
No claim is made that an arbitrary non-renewal, history-dependent reset schedule obeys the same weighted-average reduction.

\section{Supplementary derivations and finite checks}

\subsection*{B.1 Generating functions}
For the symmetric path counts, let
\[
C_+(x)=\sum_{n\ge1}N_+(n)x^n,
\qquad
C_-(x)=\sum_{n\ge1}N_-(n)x^n.
\]
The second-order path recurrence yields
\[
C_+(x)=\frac{U_{z-1}(1/(2x))}{U_{2m-1}(1/(2x))},
\qquad
C_-(x)=\frac{U_{2m-z-1}(1/(2x))}{U_{2m-1}(1/(2x))},
\]
where $U_k$ is the Chebyshev polynomial of the second kind.  Using
\[
\sin((2m-z)\theta)-\sin(z\theta)
=2\cos(m\theta)\sin((m-z)\theta)
\]
gives
\[
C_-(x)-C_+(x)=
\frac{U_{m-z-1}(1/(2x))}{U_{m-1}(1/(2x))}.
\]
This is exactly the generating function for lower first-passage paths that avoid the midpoint.  This supplies a separate generating-function description of the paths that do and do not visit the midpoint; no monotonicity statement about that conditional fraction is needed here.

\subsection*{B.2 Symmetric-walk covariance}
Let $h_z=\Pp_z(J=0)=(a-z)/a$ and $t_z=\Ee_zT=z(a-z)$ for the simple symmetric walk.  Put
\[
g_z=\Ee_z[T\ind_{\{J=0\}}].
\]
Conditioning on the first step gives
\[
g_z=\frac12(g_{z-1}+g_{z+1})
+\frac12(h_{z-1}+h_{z+1}),
\qquad g_0=g_a=0.
\]
Since $h$ is harmonic in the interior, this is
\[
g_{z+1}-2g_z+g_{z-1}=-2h_z.
\]
Solving the cubic second-difference equation and subtracting $t_zh_z$ gives
\[
\Cov_z(\ind_{\{J=0\}},T)
=-\frac{z(a-z)(a-2z)}{3a},
\]
which proves \eqref{eq:discretesusc}.

\subsection*{B.3 Exact audit algorithm}
The regression certificate propagates integer interior path counts.  At every admissible $(a,z,n)$ it checks
\[
\Delta_{a,z,n}=N_+(n+2)N_-(n)-N_+(n)N_-(n+2).
\]
All arithmetic is over Python integers; no floating-point comparison enters the discrete audit.  The underlying integer-recurrence audit can be reproduced independently from the displayed formula; the numerical audit counts are supplementary rather than hypotheses of the proofs.

\subsection*{B.4 Exact global switch-count certificate}\label[appendix]{app:exacttwoswitch}
We give a finite certificate for \cref{cor:exacttwoswitch}; all decimal endpoints below represent \emph{exact rational numbers}, not rounded floating-point results.  On the odd parity class $(1,3,5,7,9,11)$, the two-step killed matrix of the birth--death chain in \cref{prop:singlefail} is
\[
\mathsf Q=\begin{pmatrix}
\frac{12}{25}&\frac3{25}&0&0&0&0\\
\frac2{125}&\frac{2411}{2500}&\frac{49}{2500}&0&0&0\\
0&\frac{49}{100}&\frac{13}{50}&\frac14&0&0\\
0&0&\frac1{40}&\frac12&\frac{19}{40}&0\\
0&0&0&\frac3{10}&\frac{33}{50}&\frac1{25}\\
0&0&0&0&\frac9{20}&\frac1{20}
\end{pmatrix}.
\]
Its characteristic polynomial is $f(u)=u\,g(u)/200000000$, where
\[
g(u)=200000000u^5-582880000u^4+622101200u^3
-298099344u^2+63766357u-4869591.
\]
The quintic $g$ has exactly one root $\mu_j$ in each of the following disjoint intervals, ordered increasingly:
\[
\begin{array}{c|cc}
 j&\text{lower endpoint}&\text{upper endpoint}\\\hline
1&0.179111138731&0.179111138732\\
2&0.296352136606&0.296352136608\\
3&0.476408434054&0.476408434056\\
4&0.974694162199&0.974694162201\\
5&0.987834128407&0.987834128409
\end{array}
\]
Indeed, exact rational evaluation gives opposite signs of $g$ at the endpoints of each interval; since $g$ has degree five, these account for all its roots.  The additional sixth eigenvalue is zero.

Put $G(k)=(\mathsf Q^k)_{1,3}$ and $H(k)=(\mathsf Q^k)_{6,3}$, with one-based indices, and $V(k)=H(k+1)G(k)-H(k)G(k+1)$.  For $k\ge1$, write
\[
G(k)=\sum_{j=1}^{5}A_j\mu_j^k,\qquad
H(k)=\sum_{j=1}^{5}B_j\mu_j^k,
\]
where the cofactor-residue formula gives
\[
\begin{split}
A_j&=\frac{147(8000\mu_j^3-9680\mu_j^2+1820\mu_j-3)}{500000000\,f'(\mu_j)},\\
B_j&=\frac{27(62500\mu_j^2-90275\mu_j+28812)}{500000000\,f'(\mu_j)}.
\end{split}
\]
Thus
\[
V(k)=\sum_{1\le i<j\le5}c_{ij}(\mu_i\mu_j)^k,\qquad
c_{ij}=(\mu_i-\mu_j)(B_iA_j-B_jA_i).
\]
Evaluation with outward rational interval operations using the displayed root brackets proves $c_{45}>0$ and, for every other pair, $0<\mu_i\mu_j<\mu_4\mu_5$, together with the \emph{exact rational} inequality
\begin{equation}\label{eq:certifiedtail}
0\le
\sum_{(i,j)\ne(4,5)}
\frac{\overline{|c_{ij}|}}{\underline{c_{45}}}
\left(\frac{\overline{\mu_i\mu_j}}{\underline{\mu_4\mu_5}}\right)^9
<\frac{581}{1000}<1.
\end{equation}
Bars and underlines denote upper and lower rational interval bounds.  Every base ratio is less than one, so \eqref{eq:certifiedtail} holds with $9$ replaced by any integer $k\ge9$.  Hence $V(k)>0$ for every $k\ge9$.  Direct rational multiplication gives signs $(+,+,-,-,+,+)$ for $k=3,4,5,6,7,8$.  Since the first two-channel time is $n=7=2\cdot3+1$, this proves the claimed exact global count of two reversals.  The companion \texttt{verify\_paperA\_exact\_two\_switches.py} independently verifies the transition matrix, characteristic and cofactor identities, five root brackets, initial signs and strict tail inequality using Python's standard-library rational arithmetic only.

\section{Additional computational checks}

The following checks complement Section~B.3 (exact integer or rational arithmetic unless stated).
\begin{enumerate}
\item \Cref{thm:discreteMain}: independent rerun for $a\le40$, $n\le200$, $18{,}285$ comparisons, none negative; joint $\TP$ of lower-exit counts for all same-parity pairs of starting sites, $a\le24$, no violation.  Equality after opening only for $(a,z)=(4,1)$, $a\le30$.
\item \Cref{thm:odd}(a), usual order: $a\le25$ odd, eight values of $\rho\in[1/50,1/2]$, total masses from the Chebyshev generating functions, $250{,}224$ comparisons, no violation.  $u_k$: $a\le41$, $26{,}165$ comparisons; $v_k$: $a\le39$, $21{,}805$ comparisons; no violation.
\item \Cref{thm:odd}(b): $a\le39$, $p\in\{1/4,1/2,3/5,9/10\}$, $N<120$, no violation.
\item \Cref{cor:range}: $a=6$, $z=2$, $p=0.6$, $3{,}000$ random clocks; range observed $[0.3913,\,0.9998]$ against the reset-free value $0.3910$.
\item \Cref{thm:switch,prop:node,prop:singlefail}: the monotone-in-space switch and asymptotic node were rechecked numerically; the former single-turn conjecture is superseded by the exact rational counterexample in \cref{prop:singlefail}.  \Cref{prop:switchbound}: in the example of \cref{prop:singlefail}, an independent high-precision symmetrised-matrix replay gives coefficient sign sequence $+-+++---+-$ (variation $5$); exact rational recurrence gives two sign reversals, at $n=11$ and $n=15$, over $253$ successive two-step comparisons, and the spectral and exact signs agree throughout that window.  The separate rational-interval certificate in \cref{app:exacttwoswitch} now proves exactly two reversals on the infinite two-channel time domain.
\end{enumerate}

\section*{A note on the use of AI}

The authors used generative artificial intelligence systems, including OpenAI's ChatGPT (GPT-6, 4--31 October 2026) and Anthropic's Claude (Claude Opus 5.5, 5--31 October 2026), during the mathematical investigation and computational verification of this manuscript. The systems were asked to explore candidate mathematical relationships, examine stochastic-order and spectral arguments, suggest possible proof approaches, investigate reversal constructions, and assist with symbolic calculations. The tools were also used for checking mathematical exposition. In particular, the second author used Claude for a detailed mathematical reading of the arguments, for developing the spectral bound on the number of temporal reversals (\cref{prop:switchbound}), for writing the C program that generates \cref{fig:plucker}, and for checking the bibliography and the English of the text.

The systems were used for generating alternative derivations, exploring examples, debugging implementations and identifying calculations suitable for exact or independent verification. They were less reliable as authorities on the validity or generality of mathematical assertions. Suggestions and provisional conclusions were therefore treated as hypotheses requiring mathematical justification, rather than accepted as proofs. Where appropriate, the authors independently checked identities, tested counterexamples and compared computational outputs against explicit mathematical formulations.

The mathematical research, proof-search procedures and certificate criteria were developed and directed by the authors. Computational assistance was evaluated through reproducible symbolic, rational, numerical or analytical checks appropriate to each result. The authors reviewed the arguments and retained responsibility for the correctness, originality, interpretation and presentation of the manuscript. No AI-generated assertion was accepted solely on the authority of a language model.

\paragraph{Code availability.} The replay scripts used for the finite verifications are available from the authors on request and will be deposited in a public repository upon publication.


\begin{thebibliography}{99}

\bibitem{EvansMajumdar2011}
M.~R. Evans and S.~N. Majumdar,
``Diffusion with stochastic resetting,''
\emph{Phys. Rev. Lett.} \textbf{106}, 160601 (2011).
\href{https://doi.org/10.1103/PhysRevLett.106.160601}{doi:10.1103/PhysRevLett.106.160601}.

\bibitem{PalPrasad2019}
A. Pal and V.~V. Prasad,
``First passage under stochastic resetting in an interval,''
\emph{Phys. Rev. E} \textbf{99}, 032123 (2019).
\href{https://doi.org/10.1103/PhysRevE.99.032123}{doi:10.1103/PhysRevE.99.032123}.

\bibitem{VillarroelMonteroVega2021}
J. Villarroel, M. Montero, and J.~A. Vega,
``A semi-deterministic random walk with resetting,''
\emph{Entropy} \textbf{23}, 825 (2021).
\href{https://doi.org/10.3390/e23070825}{doi:10.3390/e23070825}.

\bibitem{Vega2026Critical}
J.~A. Vega Coso,
``Channel duality and midpoint invariance in random walks with resetting,''
arXiv:2603.24803 (2026).
\href{https://doi.org/10.48550/arXiv.2603.24803}{doi:10.48550/arXiv.2603.24803}.

\bibitem{Vega2026Duality}
J.~A. Vega Coso,
``Spectral duality and reset-neutral distributions in random walks with multi-site geometric resetting,''
arXiv:2605.00657 (2026).
\href{https://doi.org/10.48550/arXiv.2605.00657}{doi:10.48550/arXiv.2605.00657}.

\bibitem{Vega2026Resolvent}
J.~A. Vega Coso,
``Resolvent intertwining and spectral duality in Markov chains with geometric resetting,''
arXiv:2608.15140 (2026).
\href{https://doi.org/10.48550/arXiv.2608.15140}{doi:10.48550/arXiv.2608.15140}.

\bibitem{Vega2026Algebra}
J.~A. Vega Coso,
``Algebra of spectral duality structures,''
arXiv:2609.04430 (2026).
\href{https://doi.org/10.48550/arXiv.2609.04430}{doi:10.48550/arXiv.2609.04430}.

\bibitem{Karlin1964}
S. Karlin,
``Total positivity, absorption probabilities and applications,''
\emph{Trans. Amer. Math. Soc.} \textbf{111}, 33--107 (1964).
\href{https://doi.org/10.1090/S0002-9947-1964-0168010-2}{doi:10.1090/S0002-9947-1964-0168010-2}.

\bibitem{KarlinMcGregor1960}
S. Karlin and J. McGregor,
``Classical diffusion processes and total positivity,''
\emph{J. Math. Anal. Appl.} \textbf{1}, 163--183 (1960).
\href{https://doi.org/10.1016/0022-247X(60)90020-2}{doi:10.1016/0022-247X(60)90020-2}.

\bibitem{Belan2018}
S. Belan,
``Restart could optimize the probability of success in a Bernoulli trial,''
\emph{Phys. Rev. Lett.} \textbf{120}, 080601 (2018).
\href{https://doi.org/10.1103/PhysRevLett.120.080601}{doi:10.1103/PhysRevLett.120.080601}.

\bibitem{EvansMajumdarSchehr2020}
M.~R. Evans, S.~N. Majumdar, and G. Schehr,
``Stochastic resetting and applications,''
\emph{J. Phys. A: Math. Theor.} \textbf{53}, 193001 (2020).
\href{https://doi.org/10.1088/1751-8121/ab7cfe}{doi:10.1088/1751-8121/ab7cfe}.

\bibitem{Vidmar2022}
M. Vidmar,
``On laws exhibiting universal ordering under stochastic restart,''
\emph{Commun. Stat. Theory Methods} \textbf{51}(5), 1290--1305 (2022).
\href{https://doi.org/10.1080/03610926.2020.1759640}{doi:10.1080/03610926.2020.1759640}.

\bibitem{PalEtAl2025}
S. Pal, L. Dagdug, D. Ghosh, D. Boyer, and A. Pal,
``Universal criterion for selective outcomes under stochastic resetting,''
\emph{Phys. Rev. E} \textbf{112}, 034116 (2025).
\href{https://doi.org/10.1103/p3yc-kmt1}{doi:10.1103/p3yc-kmt1}.

\bibitem{MullerStoyan2002}
A. M\"uller and D. Stoyan,
\emph{Comparison Methods for Stochastic Models and Risks},
Wiley, Chichester (2002).

\bibitem{BartoszewiczSkolimowska2006}
J. Bartoszewicz and M. Skolimowska,
``Preservation of classes of life distributions and stochastic orders under weighting,''
\emph{Stat. Probab. Lett.} \textbf{76}(6), 587--596 (2006).
\href{https://doi.org/10.1016/j.spl.2005.09.003}{doi:10.1016/j.spl.2005.09.003}.

\bibitem{Lindstrom1973}
B. Lindstr\"om,
``On the vector representations of induced matroids,''
\emph{Bull. London Math. Soc.} \textbf{5}, 85--90 (1973).
\href{https://doi.org/10.1112/blms/5.1.85}{doi:10.1112/blms/5.1.85}.

\bibitem{GesselViennot1985}
I.~M. Gessel and G. Viennot,
``Binomial determinants, paths, and hook length formulae,''
\emph{Adv. Math.} \textbf{58}(3), 300--321 (1985).
\href{https://doi.org/10.1016/0001-8708(85)90121-5}{doi:10.1016/0001-8708(85)90121-5}.

\bibitem{PolyaSzego}
G. P\'olya and G. Szeg\H{o},
\emph{Problems and Theorems in Analysis II},
Classics in Mathematics, Springer, Berlin (1998); reprint of the 1976 edition.
\href{https://doi.org/10.1007/978-3-642-61905-2}{doi:10.1007/978-3-642-61905-2}.

\bibitem{ShakedShanthikumar2007}
M. Shaked and J.~G. Shanthikumar,
\emph{Stochastic Orders},
Springer, New York (2007).
\href{https://doi.org/10.1007/978-0-387-34675-5}{doi:10.1007/978-0-387-34675-5}.

\bibitem{GantmacherKrein}
F.~R. Gantmacher and M.~G. Krein,
\emph{Oscillation Matrices and Kernels and Small Vibrations of Mechanical Systems},
rev. ed., AMS Chelsea, Providence (2002).
\href{https://doi.org/10.1090/chel/345}{doi:10.1090/chel/345}.

\bibitem{Karlin1968}
S. Karlin,
\emph{Total Positivity}, Vol.~I,
Stanford University Press, Stanford (1968).

\bibitem{NewtonReversals}
T. Newton,
``A sharp even-interval threshold for two temporal reversals: Schur condensation and the final length-ten classification,''
working manuscript (2026).

\bibitem{PaperB}
T. Newton and J.~A. Vega Coso,
``Exit pairs as Weyl data: reset neutrality, isospectrality and rigidity,''
in preparation (2026).

\end{thebibliography}
\end{document}